\documentclass[12pt]{article}
\usepackage{amsmath,amsthm,amstext,amssymb,a4,verbatim}

\title{Conformal Geometry of
Surfaces in $S^4$\\ and Quaternions}
\author{F. Burstall, D. Ferus, K. Leschke, F. Pedit, U. Pinkall}
\date{}
\newsavebox{\address}
\sbox{\address}{%
\parbox{14cm}{\mbox{}\hspace{-3pt}
\parbox[t]{6cm}{\small{
Francis E. Burstall\\
School of Mathematical Sciences\\
University of Bath\\
GB-BATH BA2 7AY
}}
\parbox[t]{7cm}{\small{
Dirk Ferus, Katrin Leschke, Ulrich Pinkall\\
Fachbereich Mathematik\\ Technische Universit\"at Berlin\\
Str. des 17.Juni 135\\
D-10623 Berlin
}}
\\[2ex]
\mbox{}\parbox[t]{6cm}{\small{
Franz Pedit\\
Department of Mathematics\\
University of Massachusetts\\
Amherst, MA 01003, USA
}}
\parbox[t]{7cm}{\small{
f.e.burstall@maths.bath.ac.uk\\
ferus@math.tu-berlin.de\\
leschke@math.tu-berlin.de\\
pedit@gang.umass.edu\\
pinkall@math.tu-berlin.de
}}}}

\theoremstyle{plain}
\newtheorem{theorem}{Theorem}
\newtheorem*{theorem*}{Theorem}
\newtheorem{lemma}{Lemma}
\newtheorem*{lemma*}{Lemma}
\newtheorem{proposition}{Proposition}
\newtheorem*{proposition*}{Proposition}

\theoremstyle{definition}
\newtheorem*{definition}{Definition}
\newtheorem{remark}{Remark}
\newtheorem*{remark*}{Remark}
\newtheorem*{remarks*}{Remarks}
\newtheorem{example}{Example}
\newtheorem*{example*}{Example}

\numberwithin{equation}{section}

\newcommand{\RR}{\mathbb{R}}
\newcommand{\CC}{\mathbb{C}}
\newcommand{\HH}{\mathbb{H}}

\newcommand{\NN}{\mathbb{N}}

\newcommand{\cZ}{\mathcal{Z}}
\newcommand{\cH}{\mathcal{H}}

\renewcommand{\Re}{\operatorname{Re}}
\renewcommand{\Im}{\operatorname{Im}}
\newcommand{\Hom}{\operatorname{Hom}}
\newcommand{\End}{\operatorname{End}}
\newcommand{\trace}{\operatorname{trace}}
\newcommand{\rank}{\operatorname{rank}}
\newcommand{\im}{\operatorname{image }}
\newcommand{\ord}{\operatorname{ord}}
\newcommand{\ind}{\operatorname{ind}}

\newcommand{\dbar}{\bar{\partial}}
\newcommand{\tL}{\tilde{L}}
\newcommand{\tS}{\tilde{S}}
\newcommand{\tQ}{\tilde{Q}}
\newcommand{\tA}{\tilde{A}}
\newcommand{\tJ}{\tilde{J}}
\newcommand{\tdelta}{\tilde{\delta}}
\newcommand{\hL}{\hat{L}}

\newcommand{\vect}[1]{%
\begin{pmatrix}#1\end{pmatrix}}
\newcommand{\aff}[1]{%
\begin{bmatrix}#1\end{bmatrix}}

\begin{document}
\maketitle

This is the first comprehensive introduction to the authors' recent
attempts toward a  better understanding of the global concepts behind spinor
representations of surfaces in 3-space.  The important new aspect is a
quaternionic-valued function theory,  whose "meromorphic functions" are
conformal maps
into $\HH$,  which extends the classical complex function theory on  Riemann
surfaces.
The first results along these lines were presented at the ICM 98 in
Berlin \cite{icm}. Basic constructions of complex Riemann surface
theory, such as holomorphic line bundles, holomorphic curves  in projective
space, Kodaira embedding, and Riemann-Roch, carry over to the
quaternionic setting. Additionally, an important new invariant of the
quaternionic holomorphic theory is the Willmore energy. For quaternionic
holomorphic curves in $\HH P^1$ this energy is the classical Willmore energy of
conformal surfaces.

The present paper is based on a course given by one of the
authors at the Summer School on Differential Geometry at Coimbra in
September, 1999. It centers on Willmore surfaces in the conformal 4-sphere
$\HH P^1$. The first three sections introduce linear algebra over the
quaternions
and the quaternionic projective line as a model for the conformal 4-sphere.
Conformal surfaces $f:M\to\HH P^1$ are identified with the pull-back of
the tautological bundle. They are treated as quaternionic line subbundles of
the trivial bundle $M\times\HH^2$. A central object, explained in section
5,  is the
mean curvature sphere (or conformal Gauss map) of such a surface, which is a
complex structure on $M\times\HH^2$. It leads to the definition of the Willmore
energy, the critical points of which are called Willmore surfaces. In section 7
we identify the new notions of our quaternionic theory with notions in
classical
submanifold theory. The rest of the paper is devoted to applications: We
classify
super-conformal immersions as twistor projections from $\CC P^3$ in the
sense of Penrose, we construct B\"acklund transformations for Willmore
surfaces in
$\HH P^1$, we set up a duality between Willmore surfaces in $S^3$ and
certain minimal
surfaces in hyperbolic 3-space, and we give a new proof of a recent
classification
result by Montiel on Willmore 2-spheres in the 4-sphere.

\tableofcontents

\vfill

{\em Research supported by Sfb288 at TU Berlin. Franz Pedit was also
supported by NSF
grants DMS-9011083 and DMS-9705479. We thank the organizers of the Summer
School on
Differential Geometry 1999 at Coimbra, Portugal, for the opportunity to
present a
course with this material.}

\begin{section}{%
Quaternions}

\begin{subsection}{The Quaternions}\label{subsection:quaternions}
The Hamiltonian quaternions $\HH$ are
the unitary $\RR$-algebra generated by the symbols $i,j,k$  with the relations
\begin{gather*}
i^2=j^2=k^2=-1,\\ ij=-ji=k,\quad jk=-kj=i,\quad ki=-ik=j.
\end{gather*}
The multiplication is associative but obviously not commutative, and each
non-zero ele\-ment has a multiplicative inverse: We have a skew-field, and
a 4-dimensio\-nal division algebra over the reals. Frobenius showed in
1877 that $\RR,\CC$ and $\HH$ are in fact the only finite-dimensional
$\RR$-algebras that are associative and have no zero-divisors. For the
element
\begin{align}
a=a_0+a_1i+a_2j+a_3k,\quad a_l\in\RR,\label{eq:quatrep}
\end{align}
we define
\begin{align*}
\bar{a}&:=a_0-a_1i-a_2j-a_3k,\\
\Re a&:=a_0,\\
\Im a&:=a_1i+a_2j+a_3k.
\end{align*}
Note that, in contrast with the complex numbers, $\Im a$ is not a real
number, and that conjugation obeys
\begin{align*}
\overline{ab}=\bar{b}\,\bar{a}.
\end{align*}

We shall identify the real vector space $\HH$ in the obvious way with
$\RR^4$, and the subspace of purely imaginary quaternions with $\RR^3$:
\begin{align*}
\RR^3=\Im \HH.
\end{align*}
The reals are identified with $\RR1$. The embedding of the complex numbers
$\CC$ is less canonical. The quaternions
$i,j,k$ equally qualify for the complex imaginary unit, and in fact any
purely imaginary quaternion of square -1 would do the job. From now on,
however, we shall  usually use the subfield $\CC\subset\HH$ generated by
$1,i$.

Occasionally we shall need the Euclidean inner
product on
$\RR^4$ which can be written as
\begin{align*}
<a,b>_{\RR}=\Re (\bar{a}b)=\Re (a\bar{b})=\frac{1}{2}(\bar{a}b+\bar{b}a).
\end{align*}

We define
\begin{align*}
|a|:=\sqrt{<a,a>_{\RR}}=\sqrt{a\bar{a}}.
\end{align*}
Then
\begin{align}
|ab|=|a|\,|b|.\label{eq:multnorm}
\end{align}
A closer study of the quaternionic multiplication displays  nice
geo\-metric aspects.

We first mention that the quaternion
multiplication incorporates both  the usual vector and scalar products on
$\RR^3$. In fact, using the representation \eqref{eq:quatrep} one finds for
$a,b\in \Im\HH=\RR^3$
\begin{gather}
ab=a\times b\, -\,<a,b>_{\RR}.\label{cross}
\end{gather}
As a consequence we state

\begin{lemma}\label{lemma:twosphere}
For $a,b\in\HH$ we have
\begin{enumerate}
\item  \label{item:aa}
$ab=ba$ if and only if $\Im a$ and $\Im
b$ are linearly dependent over the reals. In particular, the reals are the
only quaternions that commute with all others.
\item   \label{item:bb}
$a^2=-1$ if and only if $|a|=1$ and $a=\Im a$.
Note that the set of all such $a$ is the usual two-sphere
\begin{align*}
S^2\subset\RR^3=\Im\HH.
\end{align*}
\end{enumerate}
\end{lemma}
\begin{proof}
Write $a=a_0+a',b=b_0+b'$, where the prime denotes the imaginary
part. Then
\begin{align*}
ab&=a_0b_0+a_0b'+a'b_0+a'b'\\
&=a_0b_0+a_0b'+a'b_0+a'\times b'-<a',b'>_{\RR}.
\end{align*}
All these products, except for the cross-product, are commutative, and
\eqref{item:aa} follows. From the same formula with $a=b$ we obtain $\Im
a^2=2a_0a'$. This vanishes if and only if $a$ is real or purely imaginary.
Together with
\eqref{eq:multnorm} we obtain \eqref{item:bb}.
\end{proof}
\end{subsection}

\begin{subsection}{The Group $S^3$}\label{subsection:sdrei}
The set of unit quaternions
\begin{align*}
S^3:=\{\mu\in\HH\,|\,|\mu|^2=1\}
\end{align*}
i.e. the 3-sphere in $\HH=\RR^4$, forms a group under multiplication. We
can also interpret it as the group of linear maps $x\mapsto\mu x$ of $\HH$
preserving the hermitian inner product
\begin{align*}
<a,b>:=\bar{a}b.
\end{align*}
This group is called the symplectic group $Sp(1)$.

We now consider the action of $S^3$ on $\HH$ given by
\begin{align*}
S^3\times\HH\to\HH,\quad (\mu,a)\mapsto \mu a\mu^{-1}.
\end{align*}
By \eqref{eq:multnorm} this action preserves the norm on $\HH=\RR^4$ and,
hence, the Euclidean scalar product. It obviously stabilizes
$\RR\subset\HH$ and,  therefore, its orthogonal
complement $\RR^3=\Im\HH$.
We get a map, in fact a representation,
\begin{align*}
\pi:S^3\to SO(3),\mu\mapsto
\mu\ldots\mu^{-1}|_{\Im\HH}.
\end{align*}
Let us compute the differential of $\pi$. For $\mu\in S^3$ and $v\in
T_{\mu}S^3=(\RR\mu)^\perp$, we get
\begin{align*}
d_\mu\pi(v)(a)=va\mu^{-1}-\mu a\mu^{-1}v\mu^{-1}
=\mu(\mu^{-1}va-a\mu^{-1}v)\mu^{-1}.
\end{align*}
Now $\mu^{-1}v$ commutes with all $a\in\Im\HH$ if and only if $v=r\mu$ for
some real $r$. But then $v=0$, because $v\perp\mu$. Hence $\pi$ is a local
diffeomorphism of
$S^3$ onto the 3-dimensional manifold $SO(3)$ of orientation preserving
orthogonal transformations of $\RR^3$. Since
$S^3$ is compact and $SO(3)$ is connected, this is a covering. And since
$\mu a\mu^{-1}=a$ for all $a\in\Im\HH$ if and only if $\mu\in\RR$, i.e. if
and only if $\mu=\pm 1$, this covering is 2:1. It is obvious  that
antipodal points of $S^3$ are mapped onto the same orthogonal
transformation, and therefore we see that
\begin{align*}
SO(3)\cong S^3/\{\mu\sim-\mu\}=\RR P^3.
\end{align*}
We have now displayed the group of unit quaternions as the universal
covering of $SO(3)$. This group is also called the spin group:
\begin{align*}
S^3=Sp(1)=Spin(3).
\end{align*}
If we identify $\HH=\CC\oplus\CC j=\CC^2$, we can add yet another
isomorphism:
\begin{align*}
S^3\cong SU(2).
\end{align*}
In fact, let $\mu=\mu_0+\mu_1j\in S^3$ with $\mu_0,\mu_1\in\CC$. Then for
$\alpha, \beta\in\RR$ we have $j(\alpha+i\beta)=(\alpha-i\beta)j$.
Therefore the $\CC$-linear map $A_{\mu}:\CC^2\to\CC^2,x\mapsto\mu x$ has
the following matrix representation with respect to the basis $1,j$ of
$\CC^2$:
\begin{alignat*}{4}
A_{\mu}1&=&\mu_0+\mu_1 j&=&1\mu_0+j\bar{\mu}_1\\
A_{\mu}j&=&-\mu_1 +\mu_0 j&=&1(-\mu_1)+j\bar{\mu}_0.
\end{alignat*}
Because of $\mu_0\bar{\mu}_0+\mu_1\bar{\mu}_1=1$, we have
\begin{align*}
\begin{pmatrix}
\mu_0&\bar{\mu}_1\\
-\mu_1&\bar{\mu}_0
\end{pmatrix}
\in SU(2).
\end{align*}
\end{subsection}
\end{section}

\begin{section}{%
Linear Algebra over the Quaternions}\label{sec:linalg}

\begin{subsection}{Linear Maps, Complex
Quaternionic Vector Spaces}\label{subsection:linalg}

Since we consider vector spaces
$V$ over the skew-field of quaternions, there are two options for the
multiplication by scalars. We choose quaternion vector spaces to be {\it
right} vector spaces, i.e. vectors are multiplied by quaternions from the
right:
\begin{align*}
V\times\HH\to V,(v,\lambda)\mapsto v\lambda.
\end{align*}
The notions of basis, dimension, subspace, and linear map work as in the usual
commutative linear algebra. The same is true for the matrix representation
of linear maps in finite dimensions. However,
there is no reasonable definition for the elementary symmetric functions
like trace and determinant: The linear map $A:\HH\to\HH,x\mapsto ix$,
has matrix $(i)$ when using $1$ as basis for $\HH$, but matrix $(-i)$ when
using the basis $j$.

If $A\in \End(V)$ is an endomorphism, $v\in
V$, and
$\lambda\in
\HH$ such that
\begin{align*}
Av=v\lambda,
\end{align*}
then for any $\mu\in\HH\backslash\{0\}$ we find
\begin{align*}
A(v\mu)=(Av)\mu=v\lambda\mu=(v\mu)(\mu^{-1}\lambda\mu).
\end{align*}
If $\lambda$ is real then the eigenspace is a quaternionic subspace.
Otherwise it is a real -- but {\em not} a quaternionic -- vector
subspace, and we obtain a whole 2-sphere of ``associated
eigenvalues'' (see Section~\ref{subsection:sdrei}). This is related to
the fact that multiplication by a quaternion (necessarily from the right)
is {\it not} an $\HH$-linear endomorphism of $V$. In fact, the space of
$\HH$-linear maps between quaternionic vector spaces  is {\it not} a
quaternionic vector space itself.

Any quaternionic vector space $V$ is of course a complex vector
space, but this structure depends on choosing an imaginary unit, as
mentioned in
section~\ref{subsection:quaternions}. We shall
instead (quite regularly) have an {\em additional} complex structure on $V$,
acting from the left, and hence commuting with the quaternionic structure.
In other words, we consider a fixed $J\in\End(V)$ such that $J^2=-I$.
Then
\begin{align*}
(x+iy)v:=vx+(Jv)y.
\end{align*}
In this case we call $(V,J)$ a {\it complex quaternionic (bi-)vector
space}. If
$(V,J)$ and $(W,J)$ are such spaces, then the quaternionic linear maps from
$V$ to
$W$ split  as a direct sum of the real vector spaces of complex linear
($AJ=JA$) and anti-linear $(AJ=-JA$) homomorphisms.
\begin{align*}
\Hom(V,W)=\Hom_+(V,W)\oplus\Hom_-(V,W)
\end{align*}
In fact, $\Hom(V,W)$ and $\Hom_\pm(V,W)$ are {\it complex} vector space
with multiplication given by
\begin{align*}
(x+iy)Av:=(Av)x+(JAv)y.
\end{align*}

The standard example of a quaternionic vector space is $\HH^n$.
An example of a {\em complex} quaternionic vector space is $\HH^2$ with
$J(a,b):=(-b,a)$.

On $V=\HH$, any complex structure is simply left-multiplication by
some $N\in\HH$ with $N^2=-1$. The following lemma describes a situation
that naturally produces such an $N$, and that will become a standard
situation for us. But, before stating that lemma, let us make a simple
observation:

\begin{remark}
On a real 2-dimensional vector space $U$ each complex structure
$J\in\End(U)$ induces an orientation $\mathcal{O}$ such that  $(x,Jx)$ is
positively oriented for any $x\neq 0$. We then call $J$  compatible with
$\mathcal{O}$.
\end{remark}

\begin{lemma}[Fundamental lemma]\label{lemma:cstructure}
\mbox{}\\[-20pt]
\begin{enumerate}
\item \label{item:aaa}
Let $U\subset\HH$ be a real subspace of dimension 2. Then there
exist $N,R\in\HH$ with the following three properties:
\begin{gather}
N^2=-1=R^2,\\
NU=U=UR,\\
U=\{x\in\HH\,|\,NxR=x\}.\label{eq:eigenU}
\end{gather}
The pair $(N,R)$ is unique up to sign. If $U$ is oriented, there is only
one such pair such that $N$ is compatible with the orientation.
\item \label{item:bbb}If
$U,N$ and $R$ are as above, and
$U\subset\Im\HH$, then
\begin{align*}
N=R,
\end{align*}
and this is a Euclidean unit normal vector of $U$ in $\Im\HH=\RR^3$.
\item \label{item:ccc}
Given $N,R\in\HH$ with $N^2=-1=R^2$, the
sets
\begin{align*}
U:=\{x\in\HH\,|\,NxR=x\},\quad U^\perp:=\{x\in\HH\,|\,NxR=-x\}
\end{align*}
are orthogonal real subspaces of dimension 2.
\end{enumerate}
\end{lemma}

\begin{definition}
Motivated by \eqref{item:bbb} of the lemma, $N$ and $R$ are called a {\em
left} and {\em  right normal vector} of $U$, though in general they are
{\em not at all} orthogonal to $U$ in the geometric sense.
\end{definition}
\begin{proof}[Proof of the lemma]
\eqref{item:aaa}. If $1\in U$ and if $a\in U$ is a unit vector orthogonal to 1,
then $a^2=-1$. Hence $(N,R)=(a,-a)$ works for $U$, and the uniqueness, up
to sign,
follows easily from $N1\in U$ and $Na\in U$. If $U$ is arbitrary, and
$x\in U\backslash\{0\}$ then put $\tilde{U}:=x^{-1}U$. Clearly,
$1\in\tilde{U}$.
Moreover, $(N,R)$ works for $U$ if and only if $(x^{-1}Nx,R)$ works for
$\tilde{U}$.

\eqref{item:bbb}.  If $U\subset\Im\HH=\RR^3$, and  $u,v$ is an orthonormal
basis of $U$, then $N=R=u\times v=uv$ satifies the requirements: Use the
geometric properties of the cross product.

\eqref{item:ccc}.  The above argument shows that $\sigma(x):=NxR$ has $\pm
1$-eigenspaces of real dimension 2. Since $\sigma$ is orthogonal, so are
its eigenspaces.
\end{proof}

\begin{example}\label{example:homplus}
 Let $(V,J),(W,J)$ be complex
quaternionic vector spaces of dimension 1. Then $\Hom_+(V,W)$ is of real
dimension 2. To see this, choose bases $v$ and $w$, and assume
\begin{align*}
Jv=vR,\quad Jw=wN.
\end{align*}
Then $N^2=-1=R^2$. Now $F\in\Hom(V,W)$ is given by $F(v)=wa$, and
\begin{align*}
FJ=JF&\iff FJv=JFv\\
&\iff  waR=J(wa)=(Jw)a=wNa \iff aR=Na.
\end{align*}
But the set of all such $a$ is of real dimension 2, by the last part of the
lemma. The same result holds for $\Hom_-(V,W)$. As stated earlier,
$\Hom_{\pm}(V,W)$ are complex vector spaces, and therefore
(non-canonically) isomorphic with $\CC$.
\end{example}
\end{subsection}

\begin{subsection}{Conformal Maps}\label{subsection:confmaps}

A linear map $F:V\to W$ between Euclidean vector spaces is called {\it
conformal} if there exists a positive $\lambda$ such that
\begin{align*}
<Fx,Fy>=\lambda<x,y>
\end{align*}
for all $x,y\in V$. This is equivalent to the fact that $F$ maps a
normalized orthogonal basis of $V$ into a normalized orthogonal basis of
$F(V)\subset W$. Here ``normalized'' means that all vectors have the same
length, possibly $\neq1$.

If $V=W=\RR^2=\CC$, and $J:\CC\to\CC$ denotes multiplication by the
imaginary unit, then $J$ is orthogonal. For $x\in\CC,|x|\neq0,$ the
vectors $(x,Jx)$ form a normalized orthogonal basis. The map $F$ is
conformal if and only if $(Fx,FJx)$
is again normalized orthogonal. On the other hand
$(Fx,JFx)$ {\em is} normalized orthogonal. Hence $F$ is conformal, if and
only if
\begin{align*}
FJ=\pm JF,
\end{align*}
where the sign depends on the orientation behaviour of $F$.

Note that this condition does not involve the scalar product, but only
involves the
complex structure $J$. A generalization of this fact to quaternions is
fundamental for the theory presented here.

If  $F:\RR^2=\CC \to\RR^4=\HH$  is $\RR$-linear and injective, then
$U=F(\RR^2)$ is a real 2-dimensional subspace of $\HH$, oriented by $J$.
Let $N,R\in\HH$ be its left and right normal vectors. Then
$NU=U=UR$, and $N$ induces an orthogonal endomorphism of $U$ compatible
with the Euclidean scalar product of $\RR^4$. The map
$F:\RR^2\to U$ is conformal if and only if $FJ= NF$. Hence $F:\CC\to
\HH$ is conformal if and only if there exist $N,R\in\HH,N^2=-1=R^2,$ such
that
\begin{align*}
*F:=FJ=NF=-FR.
\end{align*}

This leads to the following fundamental

\begin{definition}
Let $M$ be a Riemann surface, i.e. a 2-dimensional mani\-fold endowed with
a complex structure $J:TM\to TM,J^2=-I$. A map $f:M\to\HH=\RR^4$ is
called {\em conformal}, if there exist $N,R:M\to \HH$ such that with
$*df:=df\circ J$,
\begin{gather}
N^2=-1=R^2\label{eq:reflex}\\
\boxed{*df=Ndf=-dfR.}\label{eq:CR}
\end{gather}
If $f$ is an immersion then \eqref{eq:reflex} follows from
\eqref{eq:CR}, and $N$ and
$R$ are unique,  called the {\em left} and {\em right normal vector} of
$f$.
\end{definition}
\begin{remarks*}
\begin{itemize}
\item Equation \eqref{eq:CR} is an analog of
\begin{align*}
*df=idf
\end{align*}
for functions $f:\CC\to\CC$, i.e. of the Cauchy-Riemann equations. In
this sense conformal maps into $\HH$ are a generalization of complex
holomorphic maps.
\item If $f$ is an immersion, then $df(T_pM)\subset\HH$ is a
2-dimen\-sional real subspace. Hence, according to
Lemma~\ref{lemma:cstructure}, there exist
$N,R$, inducing a complex structure $J$ on $T_pM\cong df(T_pM)$. The
definition   requires that $J$ coincides with the complex structure
already given on $T_pM$.
\item For an immersion $f$ the existence of $N:M\to\HH$ such that $*df=Ndf$
already implies that the immersion $f:M\to\HH$ is conformal. Similarly for
$R$.
\item If $f:M\to\Im\HH=\RR^3$ is an immersion then $N=R$ is
``the classical'' unit normal vector of $f$. But for general $f:M\to\HH$,
the vectors $N$ and $R$ are {\em not} orthogonal to
$df(TM)$.
\end{itemize}
\end{remarks*}
\end{subsection}

\end{section}

\begin{section}{%
Projective Spaces}\label{sec:hp}

In complex function theory the Riemann sphere $\CC P^1$ is more convenient
as a target space for holomorphic functions than the complex plane.
Similarly, the natural target space for conformal immersions is
$\HH P^1$, rather than $\HH$. We therefore give a description of the
quaternionic projective space.

\begin{subsection}{Projective Spaces and Affine Coordinates.}
\label{subsec:aff}
The quaternionic projective space $\HH P^n$
is defined, similar to its real and complex cousins, as the set of
quaternionic lines in $\HH^{n+1}$. We have the (continuous) cano\-nical
projection
\begin{align*}
\pi:\HH^{n+1}\backslash\{0\}\to \HH P^n, x\mapsto \pi(x)=[x]=x\HH.
\end{align*}
The manifold structure of $\HH P^n$ is defined as follows:

For any linear form $\beta\in(\HH^{n+1})^*,\beta\neq0,$
\begin{align*}
u:\pi(x)\mapsto x<\beta,x>^{-1}
\end{align*}
is well-defined and maps the open set $\{\pi(x)\,|\,<\beta,x>\neq 0\}$ onto
the affine hyperplane $\beta=1$, which is isomorphic to $\HH^n$.
Coordinates of this type are called {\em affine coordinates} for $\HH P^n$.
They define a (real-analytic) atlas for $\HH P^n$.

We shall often use this in the following setting: We choose a basis for
$\HH^{n+1}$ such that $\beta$ is the last coordinate function. Then we get
\begin{align*}
\aff{x_1\\\vdots\\x_n\\x_{n+1}}\mapsto
\vect{x_1x_{n+1}^{-1}\\\vdots\\x_nx_{n+1}^{-1}\\1}
\text{ or }
\vect{x_1x_{n+1}^{-1}\\\vdots\\x_nx_{n+1}^{-1}}
\end{align*}

The set
\begin{align*}
\{\pi(x)\,|\,<\beta,x>=0\}
\end{align*}
is called the {\em hyperplane at infinity}.

\begin{example*}
In the special case $n=1$, the hyperplane at infinity is a
single point: $\HH P^1$ is the one-point compactification of $\RR^4$,
hence ``the'' 4-sphere:
\begin{align*}
\boxed{\HH P^1=S^4.}
\end{align*}
Note however, that the notion of {\em the antipodal map} is
natural on the usual 4-sphere, but not on $\HH P^1$ -- unless we introduce
additional structure, like a metric.
\end{example*}

For our purposes it is important to have a good description of the tangent
space $T_l\HH P^n$ for $l\in\HH P^n$. For that purpose, we consider the
projection
\begin{align*}
\pi:\HH^{n+1}\backslash\{0\}\to\HH P^n
\end{align*}
in affine coordinates: If
$\beta\in(\HH^{n+1})^*$ is as above, then
\begin{align*}
h=u\circ\pi:\HH^{n+1}\backslash\{0\}\to\HH^{n+1},x\to x<\beta,x>^{-1}
\end{align*}
satisfies
\begin{align*}
d_xh(v)=v<\beta,x>^{-1}-x<\beta,x>^{-1}<\beta,v><\beta,x>^{-1}.
\end{align*}
Therefore
\begin{gather*}
\ker d_xh=x\HH,\\
d_{x\lambda}h(v\lambda)=d_xh(v)
\end{gather*}
for $\lambda\in\HH\backslash\{0\}$, and the same holds for $\pi$:
\begin{gather}
\ker d_x\pi=x\HH,\label{eq:kernh}\\
d_{x\lambda}\pi(v\lambda)=d_x\pi(v).\label{eq:homo}
\end{gather}

By \eqref{eq:kernh}, $d_x\pi$
induces an isomorphism
\begin{align*}
d_x\pi:\HH^{n+1}/l\overset{\cong}{\to} T_l\HH P^n, \quad l=\pi(x),
\end{align*}
of real vector spaces, but it depends on the choice of $x\in l$. To
eliminate this dependence, we note that by \eqref{eq:homo} the map
\begin{align*}
\Hom(l,\HH^{n+1}/l)\to T_l\HH P^n, F\mapsto d_x\pi(F(x)),
\end{align*}
with $ x\in l\backslash\{0\}$ is a well-defined isomorphism:
\begin{align}
\Hom(l,\HH^{n+1}/l)\cong T_l\HH P^n.\label{eq:TPHn}
\end{align}
 In other
words, we identify $d_x\pi(v)$ with the homomorphism from $l=\pi(x)=x\HH$
to $\HH^{n+1}/l$ that maps $x$ to $\pi_l(v):=v\mod l$. For practical use,
we rephrase this as follows:
\begin{proposition}\label{proposition:deltaf}
Let  $\tilde{f}:M\to \HH^{n+1}\backslash\{0\}$ and
$f=\pi\tilde{f}:M\to\HH P^n$.\\ Let $p\in M,l:=f(p), v\in T_pM$. Then
\begin{align*}
d_pf:T_pM\to T_{f(p)}\HH P^n=\Hom(f(p),\HH^{n+1}/f(p))
\end{align*}
is given by
\begin{align*}
d_pf(v)(\tilde{f}(p)\lambda)=\pi_l(d_p\tilde{f}(v)\lambda).
\end{align*}
We denote the differential in this interpretation by
$\delta f$:
\begin{align}
\boxed{\delta f(v)(\tilde{f})=d\tilde{f}(v)\mod f.}\label{eq:deltadiff}
\end{align}
\end{proposition}
\begin{proof} The tangent vector
\begin{align*}
d_pf(v)=d_{\tilde{f}(p)}\pi(d_p\tilde{f}(v))\in T_{f(p)}\HH P^n
\end{align*}
is identified with the homomorphism $F:f(p)\to\HH^n/f(p)$, that maps
$\tilde{f}(p)$ into $d_p\tilde{f}(v)\mod f(p)$.
\end{proof}
\end{subsection}

\begin{subsection}{Metrics on $\HH P^n$.} \label{subsection:metrics}
Given a non-degenerate quaternionic her\-mi\-tian
inner pro\-duct
$<.,.>$ on $\HH^{n+1}$, we  define a (possibly degenerate
Pseudo-) Riemannian metric on $\HH P^n$ as follows: For
$x\in\HH^{n+1}$ with $<x,x>\neq0$ and
$v,w\in (x\HH)^\perp$ we define
\begin{align*}
<d_x\pi(v),d_x\pi(w)>=\frac{1}{<x,x>}\Re <v,w>.
\end{align*}
This is well-defined since, for $0\neq\lambda\in\HH$, we have
\begin{align*}
<d_{x\lambda}\pi(v\lambda),d_{x\lambda}\pi(w\lambda)>
=<d_x\pi(v),d_x\pi(w)>.
\end{align*}

It extends to arbitrary $v,w$ by
\begin{align}
<d_x\pi(v),d_x\pi(w)>
&=\Re\frac{<v,w-x<x,w><x,x>^{-1}>}{<x,x>}\nonumber\\
&=\Re\frac{<v,w><x,x>-<v,x><x,w>}{<x,x>^2}.\label{eq:met}
\end{align}

\begin{example}
For $<v,w>=\sum \bar{v_k}w_k$ we obtain the standard Rie\-mannian
metric on $\HH P^n$. (In the complex case, this is the so-called
Fubini-Study metric.) The corresponding conformal structure is in the
background
of all of the following considerations.
\item We take this standard Riemannian metric on $\HH P^1=S^4$ and ask
which metric it induces on $\RR^4$ via the  affine parameter
\begin{align*}
h:\HH\to\HH P^1,x\mapsto\aff{x\\1}.
\end{align*}
Let $\tilde{h}:\HH\to\HH^2,x\mapsto (x,1)$, and let ``$\equiv$'' denote
equality mod
$\vect{x\\1}\HH$.

Then
\begin{align*}
\delta_xh(v)(\vect{x\\1})
\equiv d_x\tilde{h}(v)
\equiv \vect{v\\0}
\equiv \vect{v\\0}-\vect{x\\1}\frac{\bar{x}v}{1+x\bar{x}}
\equiv \vect{v\\-\bar{x}v}\frac{1}{1+x\bar{x}}.
\end{align*}
The latter vector is $<.,.>$-orthogonal to $(x,1)$, and therefore the
induced metric on $\HH$ is given by
\begin{align*}
h^*<v,w>_x&=\frac{1}{(1+x\bar{x})^3}\Re<\vect{v\\-\bar{x}v},
\vect{w\\-\bar{x}w}>\\
&=\frac{1}{(1+x\bar{x})^2}\Re(\bar{v}w)
=\frac{1}{(1+x\bar{x})^2}<v,w>_{\RR}.
\end{align*}
But  stereographic
projection of $S^4$ induces the metric
\begin{align*}
\frac{2}{(1+x\bar{x})^2}<v,w>_{\RR}
\end{align*}
on $\RR^4$. Hence the standard metric on
$\HH P^1$ is of constant curvature 4.
\end{example}

\begin{example}\label{example:hyperbolic}
If we consider  an {\em indefinite} hermitian metric on $\HH^{n+1}$, then
the above construction of a metric on $\HH P^n$ fails for isotropic lines
($<l,l>=0$), but these points are scarce. We consider the case $n=1$, and the
hermitian inner product
\begin{align*}
<v,w>=\bar{v_1}w_2+\bar{v_2}w_1.
\end{align*}
Isotropic lines are characterized in affine coordinates
$h:x\mapsto\vect{x\\1}$ by
\begin{align*}
0=<\vect{x\\1},\vect{x\\1}> =\bar{x}+x,
\end{align*}
i.e. by $x\in\Im\HH=\RR^3$.

The point at infinity $\vect{1\\0}\HH$ is isotropic, too. Therefore, the set of
isotropic points  is a 3-sphere $S^3\subset S^4$, and its
complement consists of two open discs or -- in affine coordinates --
two open half-spaces.

As in the previous example, we find
\begin{align*}
h^*<v,w>_x=\frac{1}{(2\Re{x})^2}\Re(\bar{v}w)
=\frac{1}{(2\Re{x})^2}<v,w>_{\RR}
\end{align*}
for the induced metric on the
half-spaces $\Re\neq 0$ of $\HH$. This is -- up to a constant factor -- the
standard hyperbolic metric on these half-spaces.
\end{example}
\end{subsection}

\begin{subsection}{Moebius Transformations on $\HH P^1$.}
The group $Gl(2,\HH)$ acts on
$\HH P^1$ by $G(v\HH):=Gv\HH$. The kernel of this action, i.e. the set of
all $G\in Gl(2,\HH)$ such that $Gv\in v\HH$ for all $v$,  is
$\{\rho I\,|\,\rho\in\RR\}$.

How is this action compatible with the
metric induced by a positive definite hermitian metric of $\HH^2$?
Using \eqref{eq:met} we find
\begin{align*}
|dG(d_x\pi(v\lambda))|^2
&=\Re\frac{
<G(v\lambda),G(v\lambda)>-
     <G(v\lambda),Gx><Gx,G(v\lambda)>}{<Gx,Gx>^2}
\\
&=|\lambda|^2\Re\frac{
<Gv,Gv>- <Gv,Gx><Gx,Gv>}{<Gx,Gx>^2}
\\
&=|\lambda|^2|dG(d_x\pi(v))|^2
\end{align*}
Taking $G=I$ we see that for $v\neq 0$ the map
\begin{align*}
\HH\to T_{\pi(x)}\HH P^1, \lambda\mapsto d_x\pi(v\lambda)
\end{align*}
is length-preserving up to a constant factor, i.e. is a conformal
isomorphism. But the same is obviously true for the metric induced by
the pull-back under an arbitrary $G$, and therefore
$GL(2,\HH)$ acts conformally on
$\HH P^1=S^4$. We call these transformations the {\it Moebius
transformations} on $\HH P^1$. In affine coordinates they are given by
\begin{align*}
\begin{pmatrix}a&b\\c&d\end{pmatrix}\aff{x\\1}
=\aff{ax+b\\cx+d}=\aff{(ax+b)(cx+d)^{-1}\\1}.
\end{align*}
This emphasises the analogy with the complex case.

It is known that this is the full group of all orientation
preserving conformal diffeomorphisms of $S^4$, see \cite{kulkarni}.
\end{subsection}

\begin{subsection}{2-Spheres in $S^4$.}\label{subsection:twospheres}
We consider the set
\begin{align*}
\cZ=\{S\in\End(\HH^2)\,|\,S^2=-I\}.
\end{align*}
For $S\in\cZ$ we define
\begin{align*}
S':=\{l\in\HH P^1\,|\,Sl=l\}.
\end{align*}
We want to show

\begin{proposition}
\begin{enumerate}
\item $S'$ is a 2-sphere in $\HH P^1$, i.e corresponds to
a real 2-plane in $\HH=\RR^4$ under a suitable affine coordinate.
\item Each 2-sphere can be obtained in this way by an $S\in\cZ$, unique up to
sign.
\end{enumerate}
\end{proposition}
\begin{proof}
We consider $\HH^2$ as a (right) complex vector space with imaginary unit
$i$. Then $S$ is  $\CC$-linear and has a (complex) eigenvalue $N$.  If
$Sv=vN$, then
\begin{align*}
S(v\HH)=vN\,\HH=v\HH.
\end{align*}
Hence $S'\neq\emptyset$.

We choose a basis $v,w$ of $\HH^2$ such that $v\HH\in S'$, i.e. $Sv=vN$
for some $N$, and $Sw=-vH-wR$. Then $S^2=-I$ implies
\begin{align*}
N^2=-1=R^2,\quad NH=HR.
\end{align*}
For the affine parametrization $h:\HH\to\HH P^1,x\mapsto [vx+w]$ we get:
\begin{align*}
[vx+wx]\in S'
&\iff
\exists_{\gamma}\;S(vx+w)=(vx+w)\gamma\\
&\iff
\exists_{\gamma}\;vNx-vH-wR=vx\gamma+w\gamma\\
&\iff
\exists_{\gamma}\;\begin{cases}
Nx-H=x\gamma\\-R=\gamma\end{cases}\\
&\iff
Nx+xR=H.
\end{align*}

This is a real-linear equation for $x$, with associated homogeneuos equation
\begin{align*}
Nx+xR=0.
\end{align*}
By Lemma~\ref{lemma:cstructure} this is of real dimension 2, and any real
2-plane can be realized this way.
\end{proof}

Obviously, $S$ and $-S$ define the same 2-sphere. But  $S$
determines $(N,R)$, thus fixing an orientation of the above
real 2-plane and thereby of $S'$. Hence the lemma can be paraphrased
as follows:
\begin{align*}
\cZ \text{ is the set of oriented 2-spheres in }S^4=\HH P^1.
\end{align*}
\end{subsection}

\end{section}

\begin{section}{%
Vector Bundles}

We shall need vector bundles over the quaternions, and therefore  briefly
introduce
them.

\begin{subsection}{Quaternionic Vector Bundles} A quaternionic vector
bundle
$\pi:V\to M$ of rank $n$ over a smooth manifold $M$ is a real vector bundle of
rank $4n$ together with a smooth fibre-preserving action  of $\HH$ on $V$ from
the right such that the fibres become quaternionic vector spaces.

\begin{example}
The product bundle
$\pi:M\times\HH^n\to M$ with the projection on the first factor and the
obvious vector space structure on each fibre $\{x\}\times\HH^n$ is also
called {\em  the trivial bundle}.
\end{example}
\begin{example}
The points of the projective space $\HH P^n$ are the 1-dimensional
subspaces of $\HH^{n+1}$. The {\em tautological bundle}
\begin{align*}
\pi_\Sigma:\Sigma\to\HH P^n
\end{align*}
is the line bundle with
$\Sigma_l=l$. More precisely
\begin{align*}
\Sigma&:=\{(l,v)\in\HH P^n\times\HH^{n+1}\,|\,v\in l\},\\
\pi_\Sigma&:\Sigma\to \HH P^n,(l,v)\mapsto l.
\end{align*}
The differentiable and vector space structure are the obvious ones.
\end{example}

\begin{example}
If $V\to \tilde{M}$ is a quaternionic vector bundle over
$\tilde{M}$, and $f:M\to\tilde{M}$ is a map, then the ``pull-back''
$f^*V\to M$ is defined by
\begin{align*}
f^*V:=\{(x,v)\,|\,v\in V_{f(x)}\}\subset M\times V
\end{align*}
with the obvious projection and vector bundle
structure. The fibre over $x\in M$ is just the fibre of $V$ over $f(x)$.
\end{example}

We shall be concerned with maps  $f:M\to\HH P^n$ from a surface into the
projective space. To $f$ we associate the bundle
$L:=f^*\Sigma$, whose fibre  over $x$ is
$f(x)\subset\HH^{n+1}=\{x\}\times\HH^{n+1}$. The bundle $L$ is a line
subbundle of
the product bundle
\begin{align*}
H:=M\times\HH^{n+1}.
\end{align*}
Conversely, every line subbundle $L$ of $H$ over $M$ determines a map
$f:M\to\HH P^n$ by $f(x):=L_x$. We obtain an identification
\begin{center}
\boxed{\parbox{10cm}{%
     \parbox{4cm}{%
       \begin{center}Maps\\$f:M\to \HH P^n$
        \end{center}}
\hfill$\leftrightarrow$ \hfill
     \parbox{4cm}{%
       \begin{center}Line subbundles\\$L\subset H=M\times\HH^{n+1}$
        \end{center}}
                     }}
\end{center}

All natural constructions for vector spaces
extend,  fibre-wise, to operations in the category of vector bundles. For
example,
a subbundle $L$ of a vector bundle $H$ induces a quotient bundle
$H/L$ with fibres $H_x/L_x$. Given two quaternionic vector bundles
$V_1,V_2$ the {\em real} vector bundle $\Hom(V_1,V_2)$ has the fibres
$\Hom(V_{1x},V_{2x})$.  A section $\Phi\in\Gamma(\Hom(V_1,V_2))$ is called
a vector bundle homomorphism. It is a smooth map $\Phi:V_1\to V_2$ such
that for all $x$ the restriction $\Phi|_{V_{1x}}$ maps $V_{1x}$
homomorphically into $V_{2x}$. There is an obvious notion of {\em
isomorphism} for vector bundles.

\begin{example} Over $\HH P^n$ we have the product bundle
$H=\HH P^n\times\HH^{n+1}$ and, inside it, the tautological subbundle $\Sigma$.
Then
\begin{align*}
T\HH P^n\cong\Hom(\Sigma,H/\Sigma),
\end{align*}
see \eqref{eq:TPHn}.
\end{example}

\begin{example}[and Definition]\label{example:delta}
Let $L$ be a line subbundle of $H=M\times\HH^{n+1}$. Let $\pi_L:H\to
H/L\in\Gamma(\Hom(H,H/L))$ be the projection.  A section
$\psi\in\Gamma(L)\subset\Gamma(H)$ is a particular map $\psi:M\to\HH^{n+1}$. If
$X\in T_pM$, then $d\psi(X)\in H_p=\HH^{n+1}$, and
\begin{align*}
\pi_L(d\psi(X))\in (H/L)_p=\HH^{n+1}/L_p.
\end{align*}
Let $\lambda:M\to\HH$. Then
\begin{align*}
\pi_L(d(\psi\lambda)(X))=\pi_L(d\psi(X)\lambda+\psi d\lambda(X))
=\pi_L(d\psi(X))\lambda.
\end{align*}
We see that
\begin{align*}
\psi\mapsto \pi_L(d\psi(X))=:\delta(X)(\psi)
\end{align*}
is tensorial in
$\psi$, i.e. we obtain
\begin{align*}
\delta(X)=\delta_L(X)&\in \Hom(L_p,(H/L)_p).
\end{align*}
Of course this is $\RR$-linear in $X$ as well, so $\delta$ should be viewed
as a 1-form on $M$ with values in $\Hom(L,H/L)$:
\begin{align}
\boxed{\delta\in\Omega^1(\Hom(L,H/L)).}\label{eq:deltaL}
\end{align}
Let us repeat: Given $p\in M, X\in T_pM$, and $\psi_0\in L_p$, there is
a section $\psi\in\Gamma(L)$ such that $\psi(p)=\psi_0$.
Then
\begin{align*}
\boxed{\delta_p(X)\psi_0=\pi_L(d_p\psi(X))=d_p\psi(X)\mod L_p.}
\end{align*}

Note the  similarity to the second fundamental form
\begin{align*}
\alpha(X,Y)=(dY(X))^\perp.
\end{align*}
of a submanifold $M$ in Euclidan space. In the case at hand, $L$ corresponds to
$TM$ and $\HH^{n+1}/L$ corresponds to the normal bundle. This  is {\em the}
general method to measure the change of a subbundle $L$ in a (covariantly
connected) vector bundle $H$.

We can view $L$ as a map $f:M\to \HH P^n$. Even if this is an
immersion, $\delta$ clearly has  nothing to do with the second fundamental
form of $f$. Instead, comparison with Proposition~\ref{proposition:deltaf}
shows that
\begin{align*}
\delta:TM\to\Hom(L,H/L)
\end{align*}
corresponds to the derivative of
$f$, and we shall therefore call it the {\em derivative of $L$}.
\end{example}

\begin{example}
The dual $V^*:=\{\omega:V\to\HH|\omega\;
\HH\text{-linear}\}$ of a quaternionic vector space
$V$ is, in a natural way, a {\it left} $\HH$-vector space. But since we
choose quaternionic vector spaces to be {\it right} vector spaces, we use
the opposite structure: For $\omega\in V^*$ and $\lambda\in\HH$ we define
\begin{align*}
\omega\lambda:=\bar{\lambda}\omega.
\end{align*}
This  extends to quaternionic vector bundles. E.g., if $L$ is a line
bundle, i.e. of rank 1, then $L^*$ is another quaterionic line bundle,
usually denoted by $L^{-1}$.
\end{example}

A quaternionic vector bundle is called trivial if it is isomorphic with
the product bundle $M\times\HH^n$, i.e. if there exist global sections
$\phi_1,\ldots,\phi_n:M\to V$ that form a basis of the fibre everywhere.
Note that for a quaternionic {\em line} bundle over a surface the total
space $V$ has real dimension $2+4=6$, and hence any section
$\phi:M\to V$ has codimension $4$. It follows from transversality theory
that any section can be slightly deformed  so that it will not hit the
0-section. Therefore there exists a global nowhere vanishing section: Any
quaternionic line bundle over a Riemann surface is (topologically) trivial.

\end{subsection}

\begin{subsection}{Complex Quaternionic Bundles}
A {\it complex quaternionic vector bundle} is a pair $(V,J)$ consisting of
a quaternionic vector bundle $V$  and a section
$J\in\Gamma(\End(V))$ with
\begin{align*}
J^2=-I,
\end{align*}
see section \ref{subsection:linalg}.

\begin{example}
Given $f:M\to\HH, *df=Ndf$, the quaternionic line bundle $L=M\times\HH$
has a complex structure given by
\begin{align*}
Jv:=Nv.
\end{align*}
\end{example}

\begin{example}\label{example:tangS} For a given $S\in\End(\HH^2)$ with
$S^2=-I$, we identified
\begin{align*}
S'=\{l\,|\,Sl=l\}\subset\HH P^1
\end{align*}
as a 2-sphere in $\HH P^1$, see
section~\ref{subsection:twospheres}. We now compute
$\delta$, or rather the image of $\delta$, for the corresponding line bundle
$L$. In other words, we compute the tangent space of  $S'\subset\HH
P^1$.

Note that, because of $SL\subset L$, $S$ induces a complex structure on
$L$, and it also induces one (again denoted by $S$) on $H/L$ such that
$\pi_LS=S\pi_L$. Now for
$\psi\in\Gamma(L)$, we have
\begin{align*}
\delta S\psi=\pi_L d(S\psi)=\pi_LSd\psi=S\pi_Ld\psi=S\delta\psi.
\end{align*}
This shows
\begin{align*}
TS'=\im \delta\subset\Hom_+(L,H/L).
\end{align*}
But the real vector bundle
$\Hom_+(L,H/L)$ has rank 2, see Example~\ref{example:homplus}, and since
$S'$ is an embedded surface, the inclusion is an equality:
\begin{align*}
T_lS'=\Hom_+(L_l,(H/L)_l)\subset\Hom(L_l,(H/L)_l)=T_l\HH P^1.
\end{align*}
\end{example}

For our next example we generalize Lemma~\ref{lemma:cstructure}.

\begin{lemma}\label{lemma:cstr}
Let $V,W$ be 1-dimensional quaternionic vector spaces, and
\begin{align*}
U\subset\Hom(V,W)
\end{align*}
be a 2-dimensional real vector subspace. Then there exists a pair of
complex structures $J\in\End(V),\tJ\in\End(W)$, unique up to sign, such
that \begin{gather*}
\tJ U=U=UJ,\\
U=\{F\in\Hom(V,W)\,|\,\tJ F J=-F\}
\end{gather*}
If $U$ is oriented, then there is only one such pair such that
$J$ is compatible with the orientation.

Note: Here we choose the sign of $J$ in such a way that it corresponds to
$-R$ rather than $R$.
\end{lemma}
\begin{proof}
Choose non-zero basis vectors $v\in V,w\in W$. Then elements in
$\Hom(V,W)$ and endomorphisms of $V$ or  of $W$ are represented by
quaternionic $1\times1$-matrices, and therefore the assertion reduces to
that of
Lemma~\ref{lemma:cstructure}.
\end{proof}

The following is now evident:
\begin{proposition}\label{proposition:cs}
Let $L\subset H=M\times\HH^2$ be an immersed oriented surface in $\HH P^1$
with derivative $\delta\in\Omega^1(\Hom(L,H/L))$. Then there exist
unique complex structures on $L$ and $H/L$, denoted by $J,\tJ$, such that
for all $x\in M$
\begin{gather*}
\tJ\delta(T_xM)=\delta(T_xM)=\delta(T_xM)J,\\
\tJ\delta=\delta J,
\end{gather*}
and $J$  is compatible with the orientation induced by
$\delta:T_xM\to \delta(T_xM)$.
\end{proposition}
\begin{definition}
A line subbundle $L\subset H=M\times\HH^{n+1}$  over a Riemann surface $M$
is called {\em conformal} or {\em a holomorphic curve} in $\HH P^n$, if
there exists a complex structure $J$ on $L$  such that
\begin{align*}
*\delta=\delta J.
\end{align*}
\end{definition}

From the proposition we see: If $L$ is an {\em immersed} holomorphic
curve in $\HH P^1$, i.e. if $\delta$ is in addition injective, such that
$\delta(TM)\subset\Hom(L,H/L)$ is a real subbundle of rank 2, then there
 is also a complex structure $\tJ\in\Gamma(\End(H/L))$ such that
\begin{align}
*\delta=\delta J=\tJ\delta.\label{eq:holo}
\end{align}
A Riemann surface immersed into $\HH P^1$ is a holomorphic curve if and
only if the complex structures given by the proposition are compatible
with the complex structure given on $M$ in the sense of \eqref{eq:holo}.

\begin{example}
Let $f:M\to\HH$ be a conformally immersed Riemann surface with right
normal vector $R$, and let
$L$ be the line bundle corresponding to
\begin{align*}
\aff{f\\1}:M\to\HH P^1.
\end{align*}
Then
$\vect{f\\1}\in\Gamma(L)$, and
\begin{align*}
\delta (\vect{f\\1}R)&=\pi_Ld(\vect{f\\1}R)
=\pi_L(\vect{df\\0}R+\vect{f\\1}dR)\\
&=\pi_L\vect{dfR\\0}
=-\pi_L\vect{*df\\0}
=-*\delta \vect{f\\1}.
\end{align*}
If we define $J\in\End(L)$ by $J\vect{f\\1}=-\vect{f\\1}R$
then
\begin{align*}
\delta J=*\delta,
\end{align*}
hence $(L,J)$ is a holomorphic curve. Conversely, if $(L,J)$ is a
holomorphic curve, then $J\vect{f\\1}=-\vect{f\\1}R$ for some
$R:M\to\HH$, and $f$ is conformal with right normal vector $R$.
\end{example}

\end{subsection}

\begin{subsection}{Holomorphic Quaternionic Bundles}
Let $(V,J)$ be a complex quaternionic vector bundle over the Riemann
surface
$M$. We decompose
\begin{align*}
\Hom_\RR(TM,V)=KV\oplus\bar{K}V,
\end{align*}
where
\begin{align*}
KV:=\{\omega:TM\to V\,|\,*\omega=J\omega\},\\
\bar{K}V:=\{\omega:TM\to V\,|\,*\omega=-J\omega\}.\\
\end{align*}
\begin{definition}
A {\it holomorphic structure} on $(V,J)$ is a quaternionic linear map
\begin{align*}
D:\Gamma(V)\to\Gamma(\bar{K}V)
\end{align*}
such that for all $\psi\in\Gamma(V)$ and $\lambda:M\to\HH$
\begin{align}
D(\psi\lambda)=(D\psi)\lambda+\frac{1}{2}(\psi d\lambda+J\psi*d\lambda).
\label{eq:dquer}
\end{align}
A section $\psi\in\Gamma(V)$ is called {\em holomorphic} if $D\psi=0$, and
we put
\begin{align*}
H^0(V)=\ker D\subset\Gamma(V).
\end{align*}
\end{definition}

\begin{remarks*}
1. For a better understanding of this, note that for
complex-valued
$\lambda$ the anti-$\CC$-linear part (the $\bar{K}$-part) of $d\lambda$ is
given by
$\dbar\lambda=\frac{1}{2}(d\lambda+i*d\lambda)$.
In fact,
\begin{align*}
(d\lambda+i*d\lambda)(JX)
&=*d\lambda(X)-i\,d\lambda(X)=-i(d\lambda+i*d\lambda)(X).
\end{align*}
A holomorphic
structure is a generalized $\dbar$-operator. Equation
\eqref{eq:dquer} is the only natural way to make sense of a product rule
of the form ``$D(\psi\lambda)=(D\psi)\lambda+\psi\dbar\lambda$''.

2. If $L$ is a holomorphic curve in $\HH P^1$, does this mean $L$ carries
a natural holomorphic structure? This is not yet clear, but we shall come
back to this question. See also Theorem~\ref{theorem:Linverse} below.
\end{remarks*}

\begin{example} Any given $J\in\End(\HH^n), J^2=-1$, turns
$H=M\times\HH^n$ into a complex quaternionic vector bundle. Then
$\Gamma(H)=\{\psi:M\to\HH^n\}$, and
\begin{align*}
D\psi:=\frac{1}{2}(d\psi+J*d\psi)
\end{align*}
is a holomorphic structure.

\end{example}

\begin{example}\label{example:consthol}
If $L$ is a complex quaternionic line bundle and $\phi\in\Gamma(L)$ has no
zeros, then there exists exactly one holomorphic structure $D$ on $(L,J)$
such that
$\phi$ becomes holomorphic. In fact, any $\psi\in\Gamma(L)$ can be written
as
$\psi=\phi\mu$ with
$\mu:M\to\HH$, and our only chance is
\begin{align}
D\psi:=\frac{1}{2}(\phi d\mu+J\phi*d\mu).\label{eq:consthol}
\end{align}
This, indeed,  satisfies the definition of a holomorphic structure.
\end{example}
\begin{example}
If $f:M\to \HH$ is a conformal surface with left normal vector
$N$, then $N$ is a complex structure for $L=M\times\HH$, and there exists a
unique $D$ such that $D1=0$. A section $\psi=1\mu$
is holomorphic  if and only if $d\mu+N*d\mu=0$, i.e.
\begin{align*}
*d\mu=Nd\mu.
\end{align*}
The holomorphic sections are therefore the conformal maps with the same
left normal $N$ as $f$. In this case $\dim H^0(L)\ge 2$, since $1$ and $f$
are independent in $H^0(L)$.
\end{example}
\begin{theorem}\label{theorem:Linverse}
If $L\subset H=M\times\HH^{n+1}$ is a holomorphic curve with complex
structure
$J$, then the dual bundle $L^{-1}$ inherits a complex structure defined by
$J\omega:=\omega J$. The pair $(L^{-1},J)$ has a canonical holomorphic
structure
$D$ characterized by the following fact: Any quaternionic linear form
$\omega:\HH^{n+1}\to\HH$ induces a section $\omega_L\in\Gamma(L^{-1})$ by
restriction to the fibres of $L$. Then for all $\omega$
\begin{align*}
D\omega_L=0.
\end{align*}
\end{theorem}
\begin{proof} The vector bundle $L^\perp$ with fibre
$L^\perp_x=\{\omega\in(\HH^{n+1})^*\,|\;\omega|_{L_x}=0\}$ has a total
space of real dimension $4n+2$. Therefore
there exists $\omega$ such that $\omega_L$ has no zero.
Example~\ref{example:consthol} yields a unique holomorphic structure
$D$ such that $D\omega_L=0$. Now any $\alpha\in\Gamma(L^{-1})$ is of the
form
$\alpha=\omega_L\lambda$ for some $\lambda:M\to\HH$. Then, by
\eqref{eq:consthol}, for any section $\psi\in\Gamma(L)$ we have
\begin{align*}
<D&\alpha,\psi>
=\frac{1}{2}<\omega_L d\lambda+J\omega_L*d\lambda,\psi>\\
&=\frac{1}{2}(<\omega d\lambda,\psi>+<\omega*d\lambda,J\psi>)\\
&=\frac{1}{2}(d<\omega \lambda,\psi>+*d<\omega\lambda,J\psi>)
-\frac{1}{2}<\omega\lambda,d\psi+*d(J\psi)>\\
&=\frac{1}{2}(d<\alpha,\psi>+*d<\alpha,J\psi>)
-\frac{1}{2}<\omega\lambda,d\psi+*d(J\psi)>.
\end{align*}
Note that $*\delta=\delta J$ implies $d\psi+*d(J\psi)\in\Gamma(L)$, and
this allows us to replace $\omega\lambda$ by $\alpha$ in the last term as
well:
\begin{align*}
<D&\alpha,\psi>
=\frac{1}{2}(d<\alpha,\psi>+*d<\alpha,J\psi>)
-\frac{1}{2}<\alpha,d\psi+*dJ\psi>.
\end{align*}
This contains no reference to $\omega$, hence $D$ is independent of the
choice of $\omega$ such that $\omega_L$ has no zero. But the last equality
shows   $D\alpha=0$ for {\em any}
$\alpha=\omega_L$ with $\omega\in(\HH^{n+1})^*$.
\end{proof}

\begin{remark}
As we shall see in the next section, a holomorphic curve $L$ in $\HH P^1$
carries a natural holomorphic structure. In higher dimensional projective
spaces this is no longer the case. Therefore $L^{-1}$ rather than $L$ plays
a prominent role in higher codimension.
\end{remark}
\end{subsection}

\end{section}

\begin{section}{%
The Mean Curvature Sphere}

\begin{subsection}{S-Theory}\label{subsection:stheory}
Let
$M$ be a Riemann surface. Let
\begin{align*}
H:=M\times\HH^2
\end{align*}
denote the product bundle over $M$, and
let $S:M\to\End(\HH^2)\in\Gamma(\End(H))$ with $S^2=-I$
be a complex structure on $H$.
We split the differential according to type:
\begin{align*}
d\psi=d'\psi+d''\psi,
\end{align*}
where $d'$ and $d''$ denote the $\CC$-linear and anti-linear components,
respectively:
\begin{align*}
*d'=Sd',\quad *d''=-Sd''.
\end{align*}
Explicitly,
\begin{align*}
d'\psi=\frac{1}{2}(d\psi-S*d\psi),\quad
d''\psi=\frac{1}{2}(d\psi+S*d\psi).
\end{align*}
So $d''$ is a holomorphic structure on $(H,S)$, while $d'$ is an
anti-holomorphic structure, i.e. a holomorphic structure of $(H,-S)$.

In general $d(S\psi)\neq Sd\psi$, and we decompose further:
\begin{align*}
d'=\partial+A,\quad d''=\dbar+Q,
\end{align*}
where
\begin{gather*}
\partial (S\psi)=S\partial\psi,\quad\dbar (S\psi)=S\dbar\psi,\\
AS=-SA,\quad QS=-SQ.
\end{gather*}
For example, we explicitly have
\begin{align*}
\dbar\psi=\frac{1}{2}(d''\psi-Sd''(S\psi)).
\end{align*}Then $\dbar$ defines a holomorphic structure and $\partial$ an
anti-holomorphic structure on $H$, while $A$ and $Q$ are
tensorial:
\begin{align}
A&\in\Gamma(K\End_-(H)),\quad Q\in\Gamma(\bar{K}\End_-(H)).
\end{align}
For $\psi:M\to\HH^2\in\Gamma(H)$ we have, by definition of $dS$,
\begin{align*}
(dS)\psi&=d(S\psi)-Sd\psi\\
&=(\partial+A)S\psi+(\dbar+Q)S\psi
-S(\partial+A)\psi-S(\dbar+Q)\psi\\
&=AS\psi+QS\psi-SA\psi-SQ\psi\\
&=-2S(Q+A)\psi\\
&=2(*Q-*A)\psi.
\end{align*}
Hence
\begin{align}
dS=2(*Q-*A),\quad *dS=2(A-Q).\label{eq:dSQA}
\end{align}
Then
\begin{align*}
SdS=2(Q+A),
\end{align*}
whence conversely
\begin{align}
Q=\frac{1}{4}(SdS-*dS),\quad A=\frac{1}{4}(SdS+*dS).\label{eq:QA}
\end{align}
\begin{remark}
Since $A$ and $Q$ are of different type, $dS=0$ if and only if $A=0$ and
$Q=0$. If $dS=0$, then the $\pm i$-eigenspaces of the complex endomorphism
$S$ decompose $H=(M\times\CC)\oplus(M\times\CC)$. Therefore
$A$ and $Q$ measure the deviation from the "complex case".
\end{remark}
\end{subsection}

\begin{subsection}{The Mean Curvature Sphere}\label{subsection:mcs}

We now consider  an immersed holomorphic curve $L\subset
H$ in $\HH P^1$ with  derivative
$\delta=\delta_L\in\Omega^1(\Hom(L,H/L))$. Then
 there exist complex structures
$J$ on $L$ and $\tJ$ on $H/L$ such that
\begin{align*}
*\delta=\delta J=\tJ\delta.
\end{align*}

We  want to extend $J$ and $\tilde{J}$ to a complex structure of $H$, i.e.
find an
\begin{align*}
S\in\Gamma(\End(H))
\end{align*}
such that
\begin{align*}
SL=L,\quad S|_L=J,\quad \pi S=\tilde{J}\pi.
\end{align*}
Note that this implies
\begin{align*}
\pi dS(\psi)=\pi(d(S\psi)-Sd\psi)=\delta J\psi-\tilde{J}\delta\psi=0,
\end{align*}
and therefore
\begin{align}
dSL\subset L.\label{eq:dSLinL}
\end{align}

The existence of $S$ is clear: Write $H=L\oplus L'$ for some
complementary bundle $L'$. Identify $L'$ with $ H/L$ using $\pi$, and
define $S|_L:=J,S|_{L'}:=\tilde{J}$. Since $L'$ is not unique, $S$ is not
unique.  It is easy to see that
$\tilde{S}=S+R$ is another such extension if and only if
$R:M\to\End(\HH^2)$ satisfies
\begin{align*}
RH\subset L\subset \ker R,
\end{align*}
whence $R^2=0$, and
\begin{align*}
RS+SR=0.
\end{align*}
Note that $R$ can  be interpreted as an element of
$\Hom( H/L,L)$. Then $R\pi=R$.

We compute  $\tilde{Q}$:
\begin{align*}
\tilde{Q}&=\frac{1}{4}((S+R)d(S+R)-*d(S+R))\\
&=\frac{1}{4}(SdS-*dS)
+\frac{1}{4}(SdR+RdS+RdR-*dR)\\
&=Q+\frac{1}{4}(SdR+RdS+RdR-*dR).
\end{align*}
If $\psi\in\Gamma(L)$, then
\begin{gather*}
0=d(R\psi)=dR\psi+Rd\psi,\\
RdR\psi=-R^2d\psi=0
\end{gather*}
and, by \eqref{eq:dSLinL},
\begin{align*}
R\underbrace{dS\psi}_{\in\Gamma(L)}=0
\end{align*}
We can therefore continue
\begin{align*}
\tilde{Q}\psi
&=Q\psi+\frac{1}{4}(SdR\psi-*dR\psi)
 =Q\psi+\frac{1}{4}(-SRd\psi+*Rd\psi)\\
&=Q\psi+\frac{1}{4}(-SR\delta\psi+R*\delta\psi)
=Q\psi+\frac{1}{4}(-SR\delta\psi+\underbrace{R\tilde{J}}_{=RS=-SR}\delta\psi).
\end{align*}
Hence, for $\psi\in\Gamma(L)$,
\begin{align}
\tilde{Q}\psi=Q\psi-\frac{1}{2}SR\delta\psi.\label{}
\end{align}

Now we start with any extension $S$ of $(J,\tilde{J})$ and, in view of
\eqref{qschlange}, define
\begin{align}
R=-2SQ(X)\delta(X)^{-1}\pi:H\to H \label{eq:defR}
\end{align}
for some $X\neq 0$.
First note that this definition is independent
of the choice of $X\neq 0$. In fact, $X\mapsto R$  is positive-homegeneous
of degree 0, and with $c=\cos\theta,s=\sin\theta$
\begin{align*}
Q(cX+sJX)(\delta(cX+sJX))^{-1})
&=Q(X)(cI+sS)(\delta(X)(cI+sS))^{-1}\\
&=Q(X)\delta(X)^{-1}.
\end{align*}
Next
\begin{align*}
RS
&=-2SQ(X)\delta_X^{-1}\pi S
=-2SQ(X)\delta_X^{-1}\tilde{J}\pi\\
&=-2SQ(X)S\delta_X^{-1}\pi
=2S^2Q(X)\delta_X^{-1}\pi\\
&=-SR.
\end{align*}
By definition \eqref{eq:defR}
\begin{align*}
 L\subset \ker R,
\end{align*}
and from \eqref{eq:QA} and \eqref{eq:dSLinL} we get
\begin{align*}
L\supset \frac{1}{4}(SdS-*dS)L=QL,
\end{align*}
whence
\begin{align*}
RH\subset L.
\end{align*}

We have now shown that $\tilde{S}=S+R$ is another extension.
Finally, using \eqref{qschlange,} we find for $\psi\in\Gamma(L)$
\begin{align*}
\tilde{Q}\psi
&=Q\psi-\frac{1}{2}SRd\psi
=Q\psi-\frac{1}{2}S(-2SQ\delta^{-1}\pi)d\psi\\
&=Q\psi-Q\delta^{-1}\pi d\psi
=0.
\end{align*}
This shows

\begin{theorem}\label{theorem:mcs}
Let $L\subset H=M\times\HH^2$ be a holomorphic curve immersed
into
$\HH P^1$. Then there exists a unique complex structure $S$ on $H$ such
that
\begin{gather}
SL=L,\quad dSL\subset L,\label{eq:tang}\\
*\delta=\delta\circ S=S\circ\delta,\label{eq:conf}\\
Q|_L=0.\label{eq:hol}
\end{gather}
\end{theorem}

$S$ is a family of 2-spheres, a {\em sphere congruence} in
classical terms. Because $S_pL_p=L_p$ the sphere
$S_p$ goes through $L_p\in\HH P^1$, while $dSL\subset L$ (or, equivalently,
$\delta S=S\delta$) implies it is tangent to
$L$ in $p$, see examples~\ref{example:delta} and \ref{example:tangS}. In an
affine coordinate system
$\aff{f\\1}=L$ the sphere
$S_p$ has the same mean curvature vector as  $f:M\to\RR^4=\HH$ at $p$, see
Remark~\ref{remark:meancurvature}.  This motivates the
\begin{definition}
$S$ is called the {\em mean curvature sphere (congruence)} of $L$. The
diffe\-rential forms $A,Q\in\Omega^1(\End(H))$ are called  {\em the Hopf
fields} of
$L$.
\end{definition}

\begin{remark}\label{remark:holostrH}
Equations \eqref{eq:tang}, \eqref{eq:conf} imply
$d\psi+S*d\psi\in\Gamma(L)$ for $\psi\in\Gamma(L)$, whence
$d''=\dbar+Q=\frac{1}{2}(d+S*d)$ leaves
$L$ invariant. Hence an immersed holomorphic curve in $\HH P^1$ is a
holomorphic subbundle of $(H,S,d'')$ and, in particular, is a holomorphic
quaternionic vector bundle itself.
\end{remark}

\begin{example}\label{example:tangSS}
Let $S\in\End(\HH^2),S^2=-I$. Then
\begin{align*}
S'=\{l\in\HH P^1\,|\,Sl=l\}\subset\HH P^1
\end{align*}
is a 2-sphere in $\HH P^1$. Let $L$ denote the corresponding line bundle
and endow $S'$ with the complex structure inherited from the immersion.
Then the mean curvature sphere congruence of $L$ is simply the constant map
$S'\to\cZ$ of value $S$:  We have $SL=L$ by definition, and the constancy
implies  $dSL=\{0\}\subset L$ and $Q=\frac{1}{4}(SdS-*dS)=0$.
\end{example}
\end{subsection}

\begin{subsection}{Hopf Fields}\label{subsection:hopffields}

In the following we shall frequently encounter differential forms. Note that
the usual definition of the wedge product of 1-forms
\begin{align*}
\omega\wedge\theta(X,Y)=\omega(X)\theta(Y)-\omega(Y)\theta(X)
\end{align*}
can be generalized verbatim to forms $\omega_i\in\Omega^1(V_i)$ with values
in vector spaces or bundles $V_i$, provided there is a product
$V_1\times V_2\to V$. Examples are the composition
$\End(V)\times\End(V)\to\End(V)$ or the pairing between the
dual $V^*$ and $V$.

On a Riemann surface $M$, any 2-form $\sigma\in\Omega^2$ is completely
determined by the quadratic form $\sigma(X,JX)=:\sigma(X)$, and we shall,
for simplicity, often use the latter. As an example,
\begin{align*}
\omega\wedge\theta(X,JX)=\omega(X)\theta(JX)-\omega(JX)\theta(X)
\end{align*}
will be written as
\begin{align}
\omega\wedge\theta=\omega\,*\!\theta-*\omega\,\theta.
                                        \label{eq:wedgeformel}
\end{align}

We now collect some information about the Hopf fields
and the mean curvature sphere congruence $S:M\to\cZ$.

\begin{lemma}\label{lemma:dAplusQ}
\begin{align*}
d(A+Q)=2(Q\wedge Q+A\wedge A).
\end{align*}
\end{lemma}
\begin{proof}
Recall from \eqref{eq:dSQA}
\begin{align*}
SdS=2(A+Q).
\end{align*}
Therefore, using $AS=-SA,QS=-SQ$,
\begin{align*}
d(A+Q)
&=\frac{1}{2}d(SdS)
=\frac{1}{2}(dS\wedge dS)\\
&=2S(A+Q)\wedge S(A+Q)\\
&=2(A\wedge A+A\wedge Q+Q\wedge A+Q\wedge Q).
\end{align*}
But $A\wedge Q=0$ by the following {\it type argument}: Using that
\begin{align*}
\text{A is
``right $\bar{K}$'', and $Q$ ``left $\bar{K}$''}
\end{align*}
we have
\begin{align}
A\wedge Q
&=A*Q-*AQ
=A(-SQ)-(-AS)Q
=0.\label{eq:type}
\end{align}
Similarly $Q\wedge A=0$, because $A$ is left $K$ and $Q$ is right $K$.
\end{proof}

\begin{lemma}\label{lem:imageA}
Let $L\subset H $ be an immersed surface and $S$ a complex structure on $H
$ stabilizing $L$ such that $dSL\subset L$. Then $Q_{|L}=0$ is
equivalent to $AH\subset L$.
\end{lemma}
Notice that the kernels and images of the $1$-forms $A$ and $Q$ are
well-defined:
if $Q_{X}\psi=0$ for some $X\in TM$ then also $Q_{JX}\psi=-SQ_{X}\psi=0$,
and thus $Q_{Z}\psi=0$ for any $Z\in TM$. In other words, the kernels of
$Q$ and $A$ are
independent of $X\in TM$. The same remark holds for the respective images.
\begin{proof}
We first need a formula for the derivative of $1$-forms
$\omega\in\Omega^{1}(\text{End}(H))$ which stabilize $L$, i.e.,
$\omega L\subset L$. If $\pi=\pi_L$, then for $\psi\in\Gamma(L)$
\begin{align*}
\pi(d\omega(X,Y)\psi)(X,Y)
=&\pi(d(\omega\psi)(X,Y)+\omega\wedge d\psi(X,Y))\\
=&\pi
(X\cdot(\omega(Y)\psi)-Y\cdot(\omega(X)\psi)
           -\underbrace{\omega([X,Y])\psi}_{\in\Gamma(L)}\\
     &+\omega(X)d\psi(Y)-\omega(Y)d\psi(X))\\
=&\delta(X)\omega(Y)\psi
-\delta(Y)\omega(X)\psi+\pi\omega(X)d\psi(Y)-\pi\omega(Y)d\psi(X)\\
=&\delta(X)\omega(Y)\psi
-\delta(Y)\omega(X)\psi+\pi\omega(X)\delta\psi(Y)-\pi\omega(Y)\delta\psi(X)\\
=&(\delta\wedge\omega+\pi\omega\wedge\delta)(X,Y)\psi,
\end{align*}
where we wedge over composition. Note that the composition
$\pi\omega\delta$ makes sense, because $\omega(L)\subset L$, and $L$ is
annihilated by $\pi$.
We apply this to $A$ and $Q$. Since $AL\subset L,QL\subset L$ we have
\underline{on $L$}, by lemma~\ref{lemma:dAplusQ},
\begin{align*}
0
&=\frac{1}{2}\pi(Q\wedge Q+A\wedge A)
 =\pi(dA+dQ)\\
&=\delta\wedge A+\pi A\wedge\delta+\delta\wedge Q+\pi Q\wedge \delta.
\end{align*}
By a type argument similar to \eqref{eq:type}, we get $\delta\wedge A=0=\pi
Q\wedge \delta$. Further,
\begin{align*}
\pi A\wedge\delta
&=\pi A*\delta-\pi*A\delta\\
&=-2S\pi A\delta,
\end{align*}
and similarly for the remaining term. We obtain
$-\pi SA\delta=S\delta Q|_L$ or
\begin{align*}
-\pi A\delta=\delta Q|_L.
\end{align*}
Since $AL\subset L$ and $\delta(X):L\to H/L$ for $X\neq 0$ is an
isomorphism, we get $\pi A=0\iff Q|_L=0$.
\end{proof}
\end{subsection}

\begin{subsection}{The Conformal Gauss Map}
\begin{definition}
For a quaternionic vector space or bundle $V$ of rank $n$ and
$A\in\End(V)$ we define
\begin{align*}
<A>:=\frac{1}{4n}\trace_{\RR}A,
\end{align*}
where the trace is taken of the real endomorphism $A$. In particular
$<I>=1$. We obtain an indefinite scalar product $<A,B>:=<AB>$.
\end{definition}
\begin{example}
For $A=(a)$ with $a=a_0+ia_1+ja_2+ka_3\in\HH$ we have
\begin{align*}
<A>=\frac{1}{4}\,4a_0=a_0,
\end{align*}
and
\begin{align*}
<AA>=\Re a^2=a_0^2-a_1^2-a_2^2-a_3^2.
\end{align*}
\end{example}

\begin{proposition}\label{proposition:conformalgauss}
The mean curvature sphere $S$ of an immersed Riemann surface $L$ satisfies
\begin{align*}
<dS,dS>=<*dS,*dS>,\quad <dS,*dS>=0,
\end{align*}
i.e. $S:M\to\cZ$ is conformal.
\end{proposition}
Because of this proposition, $S$ is also called {\em the conformal Gauss
map}, see Bryant \cite{bryant}.
\begin{proof}
We have $QA=0$, and therefore
\begin{align}
<Q,A>=<A,Q>=0.\label{eq:qorthoa}
\end{align}
Then, from \eqref{eq:dSQA},
\begin{align*}
<dS,dS>
=&4<-S(Q+A),-S(Q+A)>
=4<Q+A,Q+A>\\
=&4<Q-A,Q-A>
=<*dS,*dS>.
\end{align*}
Similarly,
\begin{align*}
<dS,*dS>
=&4<-S(Q+A),A-Q>\\
=&4(<SQQ>-<S\underbrace{QA}_{=0}>+<SAQ>-<SAA>).
\end{align*}
But, by a property of the real trace,
\begin{align*}
<SAQ>&=<QSA>=<-SQA>=0,\\
<SQQ>&=<QSQ>=<-SQQ>=0,\\
<SAA>&=<ASA>=<-SAA>=0.
\end{align*}
\end{proof}
\end{subsection}
\end{section}

\newpage

\begin{section}{%
Willmore Surfaces}

Throughout this section $M$ denotes a {\em compact} surface.

\begin{subsection}{The Energy Functional}
The set
\begin{align*}
\cZ=\{S\in\End(\HH^2)\,|\,S^2=-I\}
\end{align*}
of oriented 2-spheres in $\HH P^1$ is a submanifold of $\End(\HH^2)$ with
\begin{align*}
T_S\cZ&=\{X\in\End(\HH^2)\,|\,XS=-SX\},\\
\perp_S\cZ&=\{Y\in\End(\HH^2)\,|\,YS=SY\}.
\end{align*}
Here we use the (indefinite) inner product
\begin{align*}
<A,B>:=<AB>=\frac{1}{8}\trace_{\RR}(AB)
\end{align*}
defined in Section~\ref{subsection:hopffields}.

\begin{definition}
The {\it energy functional} of a map $S:M\to\cZ$ of a Riemann surface $M$
is defined by
\begin{align*}
E(S):=\int_M<dS\wedge *dS>.
\end{align*}
Critical points $S$ of this functional with respect to variations of
$S$ are called {\it harmonic maps} from
$M$ to $\cZ$.
\end{definition}

\begin{proposition}\label{proposition:harmonicity}
 $S$ is harmonic if and only if the $\cZ$-tangential component of
$d*dS$ vanishes:
\begin{align}
(d*dS)^T&=0.
\end{align}
 This condition is equivalent to any of the following:
\begin{align}
d(S*dS)&=0,\\
d*A&=0,\\
d*Q&=0.
\end{align}
In fact,
\begin{align}
d(S*dS)=4d*Q=4d*A=S(d*dS)^T=(Sd*dS)^T.\label{eq:dstarQ}
\end{align}
\end{proposition}
\begin{proof}
Let $S_t$ be a variation of $S$ in $\cZ$ with
variational vector field
$\dot{S}=:Y$. Then $SY=-YS$ and
\begin{align*}
\frac{d}{dt}E(S)&=\frac{d}{dt}\int_M<dS\wedge *dS>
=\int_M<dY\wedge *dS>+<dS\wedge *dY>.
\end{align*}
Using   the wedge formula \eqref{eq:wedgeformel} and
$\trace_{\RR}(AB)=\trace_{\RR}(BA)$, we get
\begin{align*}
<dS\wedge *dY>&=<dS(-dY)-*dS*dY>=<dY\wedge*dS>.
\end{align*}
Thus
\begin{align*}
\frac{d}{dt}E(S)&=2\int_M<dY\wedge *dS>=-2\int_M<Yd*dS>=-2\int_M<Y,d*dS>.
\end{align*}
Therefore $S$ is harmonic if and only if $d*dS$ is normal.

For the other equivalences, first note
\begin{align*}
0
&=d*d(S^2)=d(*dSS+S*dS)\\
&=(d*dS)S-*dS\wedge dS+dS\wedge *dS+Sd*dS\\
&=-2(dS)^2-2(*dS)^2+(d*dS)S+Sd*dS\\
&=2dS\wedge *dS+(d*dS)S+Sd*dS.
\end{align*}
Now, together with
$*Q-*A=\frac{1}{2}dS$ and $A=\frac{1}{4}(SdS+*dS)$, this implies
\begin{align*}
8d*Q&=8d*A=2d(S*dS)\\
&=2dS\wedge*dS+2Sd*dS\\
&=-(d*dS)S+Sd*dS\\
&=S(\underbrace{d*dS+S(d*dS)S}_{=2(d*dS)^T}).
\end{align*}
\end{proof}

We now consider the case where $S$ is the mean
curvature sphere of an immersed holomorphic
curve. We decompose $dS$ into the Hopf fields.

\begin{lemma}\label{lemma:aa}
\begin{gather}
<dS\wedge *dS>=4(<A\wedge *A>+<Q\wedge *Q>),\\
<dS\wedge SdS>=4(<A\wedge *A>-<Q\wedge *Q>).
\end{gather}
\end{lemma}
\begin{proof}
Recall from section~\ref{subsection:stheory}
\begin{align*}
dS=2(*Q-*A),\quad *dS=2(A-Q),\quad SdS=2(Q+A).
\end{align*}
Further
\begin{align*}
*Q\wedge A=0,\quad *A\wedge Q=0
\end{align*}
by type. Therefore
\begin{align*}
<dS\wedge *dS>
&=4<(*Q-*A)\wedge (A-Q)>\\
&=-4<*Q\wedge Q>-4<*A\wedge A>\\
&=4<Q\wedge *Q>+4<A\wedge *A>,
\end{align*}
and similarly for $<dS\wedge SdS>$.
\end{proof}

\begin{lemma}\label{lemma:AonL}
Let $V$ be a  quaternionic vector space,
$L\subset V$ a quaternionic line, $S,B\in\End(V)$ such that
\begin{align*}
S^2=-I,\quad SB=-BS,\quad \im B\subset L.
\end{align*}
Then
\begin{align*}
\trace_{\RR}B^2\le 0,
\end{align*}
with equality if and only if $B|_L=0$.
\end{lemma}
\begin{proof} We may assume $B\neq 0$. Then $L=BV$, and
$SB=-BS$ implies $SL=L$. Let $\phi\in L\backslash\{0\}$, and
\begin{align*}
S\phi=\phi\lambda,\quad B\phi=\phi\mu.
\end{align*}
Then $\lambda^2=-1$, and $BS=-SB$ implies
\begin{align*}
\lambda\mu=-\mu\lambda.
\end{align*}
Therefore $\mu$ is imaginary, too. It follows $B^2\phi=-|\mu|^2\phi$, and
\begin{align*}
\trace_{\RR}B^2=\trace_{\RR}B^2|_L=-4|\mu|^2.
\end{align*}
\end{proof}

This can be applied to $A$ or $Q$ instead of $B$, since
their rank is $\le 1$. We obtain
\begin{lemma}\label{lemma:positive} For an
immersed holomorphic curve
$L$ we have
\begin{align}
<A\wedge*A>=\frac{1}{2}<A|_L\wedge*A|_L>, \label{eq:restrict}
\end{align}
and
\begin{align}
<A\wedge*A>\,\ge 0,\quad <Q\wedge *Q>\,\ge 0.\label{eq:posdef}
\end{align}
In particular $E(S)\ge 0$.
\end{lemma}
\begin{proof}
\begin{align*}
<A\wedge*A>&
=\frac{1}{8}\trace_{\RR}(-A^2-\underbrace{(*A)^2}_{=-ASSA=A^2})
=-\frac{1}{4}\trace_{\RR} A^2.
\end{align*}
Because  $\dim L=\frac{1}{2}\dim H$ we  similarly have
\begin{align*}
<A|_L\wedge*A|_L>=-\frac{1}{2}\trace_{\RR} A|_L^2,
\end{align*}
see section~\ref{subsection:hopffields}.
Because  $AH\subset L$, we have
\begin{align*}
\trace_{\RR}A^2=\trace_{\RR}A|_L^2.
\end{align*}
This proves \eqref{eq:restrict}.
The positivity follows from Lemma~\ref{lemma:AonL}.
\end{proof}

\begin{proposition}
\begin{enumerate}
\item The (alternating!) 2-form $\omega\in\Omega^2(\cZ)$ defined by
\begin{align*}
\omega_S(X,Y)=<X,SY>,\quad\text{for } S\in\cZ,\;X,Y\in T_S\cZ,
\end{align*}
is closed.
\item If $S:M\to\cZ$, and $dS=2(*Q-*A)$ as usual, see section
\ref{subsection:stheory} \eqref{eq:QA}, then
\begin{align*}
S^*\omega=2<A\wedge *A>-2<Q\wedge *Q>.
\end{align*}
In particular,
\begin{align*}
\deg S:=\frac{1}{\pi}\int_M <A\wedge *A>-<Q\wedge *Q>
\end{align*}
is a topological invariant of $S$.
\end{enumerate}
\end{proposition}
\begin{remark}
Since $S$ maps the surface $M$ into the 8-dimensional $\cZ$, $\deg S$
certainly is not the mapping degree of $S$. But for immersed holomorphic
curves it is the difference of two mapping degrees $\deg S=\deg N-\deg R$,
where $N,R:M\to S^2$ are the left and right normal vector in affine
coordinates, see section~\ref{section:metricgeom}.
\end{remark}
\begin{proof}
(i). We consider the 2-form on $\End(\HH^2)$ defined by
\begin{align*}
\tilde{\omega}_S(X,Y):=\frac{1}{2}(<X,SY>-<Y,SX>).
\end{align*}
Then $d_S\tilde{\omega}(X,Y,Z)$  is a linear
combination of terms of the form
\begin{align*}
<Y,XZ>.
\end{align*}
But if $X,Y,Z\in T_S\cZ,\;S\in\cZ$, we get
\begin{align*}
<Y,XZ>&=-<S^2YXZ>=<SYXZS>\\
&=<S^2YXZ>=-<Y,XZ>,
\end{align*}
hence $<Y,XZ>=0$.
Therefore, if $\iota:\cZ\to\End(\HH^2)$ is the inclusion,
\begin{align*}
d\omega=d\iota^*\tilde{\omega}=\iota^*d\tilde{\omega}=0.
\end{align*}

(ii). We have
\begin{align*}
S^*\omega(X,Y)&=<dS(X),SdS(Y)>\\
&=\frac{1}{2}(<dS(X)SdS(Y)>-<SdS(X)dS(Y)>)\\
&=\frac{1}{2}(<dS(X)SdS(Y)>-<dS(Y)SdS(X)>)\\
&=\frac{1}{2}<dS\wedge SdS>(X,Y),
\end{align*}
and Lemma~\ref{lemma:aa} yields the formula.

The topological invariance under deformations of $S$ follows from Stokes
theorem: If $\tilde{S}:M\times[0,1]\to\cZ$ deforms $S_0:M\to \cZ$ into
$S_1$, then
\begin{align*}
0&=\int_{M\times[0,1]}d\tilde{S}^*\omega\\
&=\int_{M\times1}\tilde{S}^*\omega-\int_{M\times0}\tilde{S}^*\omega\\
&=\int_MS_1^*\omega-\int_MS_0^*\omega.
\end{align*}
\end{proof}

\begin{remark}From
\begin{align*}
E(S)&=4\int_M <A\wedge *A>+<Q\wedge *Q>\\
&=8\int_M <A\wedge *A>
+\underbrace{4\int_M (<Q\wedge *Q>-<A\wedge *A>)}_{%
\text{topological invariant}}
\end{align*}
we see that for variational problems the energy functional can be replaced
by the  integral of $<A\wedge*A>$.
\end{remark}
\end{subsection}

\begin{subsection}{The Willmore Functional}
\begin{definition}
Let $L$ be a  compact immersed holomorphic curve in $\HH P^1$ with Hopf
field
$A$. The {\it Willmore functional} of $L$ is defined as
\begin{align*}
W(L):=\frac{1}{\pi}\int_M<A\wedge *A>.
\end{align*}
If we vary the immersion $L:M\to\HH P^1$, it will in general not remain
a holomorphic curve. On the other hand, any immersion {\em induces} a
complex structure $J$ on $M$ such that with respect to this it is a
holomorphic curve, see Proposition~\ref{proposition:cs}. Critical points
of $W$ with respect to such variations are called  {\it
Willmore surfaces}. If we consider only variations of $L$ fixing the
conformal structure of $M$ they are called {\em constrained Willmore
surfaces}, but we shall not treat this case here.
\end{definition}

\begin{example} For immersed surfaces in $\RR^4$ we have
\begin{align*}
W(L)=\frac{1}{4\pi}\int_M(H^2-K-K^\perp)|df|^2,
\end{align*}
see
section~\ref{subsection:willmoreaffine},
Proposition~\ref{proposition:affinwillmoreint}.
\end{example}

\begin{theorem}[Ejiri \cite{ejiri}, Rigoli \cite{rigoli}]
An immersed holomorphic curve $L$ is  Willmore  if and only if
its mean curvature sphere $S$ is harmonic.
\end{theorem}
\begin{proof}
Let $L_t$ be a variation, and $S_t$ its mean
curvature sphere. Note that for $L_t$ to stay conformal the complex
structure, i.e. the operator~$*$, varies, too.  The variation has a
variational vector field $Y\in\Gamma(\Hom(L,H/L))$ given by
\begin{align*}
Y\psi:=\pi(\left.\frac{d}{dt}\right|_{t=0}\psi),\quad\psi_t\in\Gamma(L_t).
\end{align*}
As usual, we abbreviate $\frac{d}{dt}|_{t=0}$ by a dot. Note that for
$\psi\in\Gamma(L)$
\begin{align}
\pi\dot{S}\psi&=\pi(S\psi)\dot{}-\pi S\dot{\psi}
=YS\psi-S\pi\dot{\psi}
=(YS-SY)\psi.\label{eq:sdot}
\end{align}
We now compute the variation of the energy, which is as good as the
Willmore functional as long as we vary $L$. By contrast, in the proof of
Proposition~\ref{proposition:harmonicity} the conformal structure on $M$
was fixed, and no $L$ was involved.
\begin{align*}
\left.\frac{d}{dt}\right|_{t=0} E(S_t)
&=\left.\frac{d}{dt}\right|_{t=0}\int_M <dS_t\wedge*_tdS_t>\\
&=\underbrace{\int_M <d\dot{S}\wedge*dS>}_{I}
+\underbrace{\int_M <dS\wedge\dot{*}dS>}_{II}
+\underbrace{\int_M <dS\wedge*d\dot{S}>}_{III}.
\end{align*}
In general $<A\wedge*B>=<B\wedge*A>$, because
$\trace_{\RR}(AB)=\trace_{\RR}(BA)$. Hence
\begin{align}
III=I. \label{eq:I}
\end{align}
 Next we claim
\begin{align}
II=0.\label{eq:II}
\end{align}
On $TM$ let $\dot{J}=B$, i.e. $\dot{*}\omega(X)=:\omega(BX)$. Then we
have
$BJ+JB=0$, and
\begin{align*}
<dS\wedge\dot{*}dS>(X,JX)
&=<dS(X)\dot{*}dS(JX)>-<dS(JX)\dot{*}dS(X)>\\
&=<dS(X)dS(BJX)>-<dS(JX)dS(BX)>\\
&=-<dS(X)dS(JBX)>-<dS(BX)dS(JX)>.
\end{align*}
But $S$ is conformal, see
Proposition~\ref{proposition:conformalgauss}, therefore
\begin{align*}
<dS(X)dS(JX)>=0\text{ for all }X.
\end{align*}
Differentiation with respect to $X$ yields
\begin{align*}
<dS(X)dS(JY)>+<dS(Y)dS(JX)>=0
\end{align*}
for all $X,Y$. Using this with $Y=BX$ we get \eqref{eq:II}.

Now, we compute the integral $I$.
\begin{align*}
I&=-\int_M <\dot{S},d*dS>\\
&\underset{\eqref{eq:dstarQ}}{=}4\int_M <\dot{S},Sd*Q>\\
&=\frac{1}{2}\int_M \trace_{\RR}(\dot{S}Sd*Q).
\end{align*}
We shall show in the following lemma that
\begin{align*}
\im d*Q\subset L\subset\ker d*Q.
\end{align*}
Therefore we can consider $d*Q$ as a 2-form
\begin{align*}
d*Q\in\Omega^2(\Hom(H/L,L),
\end{align*}
and continue
\begin{align*}
I&=\frac{1}{2}\int_M \trace_{\RR}(\dot{S}Sd*Q:H\to H)\\
&=\frac{1}{2}\int_M \trace_{\RR}(\pi\dot{S}Sd*Q:H/L\to H/L)\\
&=\frac{1}{2}\int_M \trace_{\RR}(\pi\dot{S}|_L Sd*Q:H/L\to H/L)\\
&\underset{\eqref{eq:sdot}}{=}
\frac{1}{2}\int_M \trace_{\RR}((YS-SY)(Sd*Q):H/L\to H/L)\\
&=-\frac{1}{2}\int_M \trace_{\RR}(Yd*Q)
  -\frac{1}{2}\int_M \trace_{\RR}(SYSd*Q).
\end{align*}
Now $d*Q$ is tangential by \eqref{eq:dstarQ}, and hence anti-commutes with
$S$. Thus
\begin{align*}
I&=
  -\frac{1}{2}\int_M \trace_{\RR}(Yd*Q)
  +\frac{1}{2}\int_M \trace_{\RR}(SYd*QS)\\
&=-\int_M \trace_{\RR}(Yd*Q)\\
&=-8\int_M <Y,d*Q>
\end{align*}
We therefore showed
\begin{align*}
\left.\frac{d}{dt}\right|_{t=0}E(S_t)=-8\int_M <Y,d*Q>.
\end{align*}
Since $\Omega^2(\Hom(H/L,L)$, this vanishes for all variational
vector fields
$Y$ if and only if
\begin{align*}
d*Q=0.
\end{align*}
\end{proof}

In the proof we made use of the following

\begin{lemma}\label{lemma:dstarQonL}
\begin{align*}
\im d*Q\subset L\subset\ker d*Q.
\end{align*}
\end{lemma}
\begin{proof}
For $\psi\in\Gamma(L)$
\begin{align*}
0=d(*Q\psi)=(d*Q)\psi-*Q\wedge d\psi=(d*Q)\psi-*Q\wedge\delta\psi,
\end{align*}
because $Q|_L=0$. But $*Q$ is right $K$, and $\delta$ is left $K$. Hence,
by type,
\begin{align*}
(d*Q)\psi=*Q\wedge\delta\psi=0.
\end{align*}
This shows the right hand inclusion.
Also,
\begin{align*}
\pi(d*Q)(X,JX)
&=\pi(d*A)(X,JX)\\
&=\pi(X\cdot(*A(JX))-JX\cdot(*A(X))
-\underbrace{*A([X,JX])}_{L-\text{valued}})\\
&=\delta(X)*A(JX)-\delta(JX)*A(X)\\
&=-\delta(X)A(X)-\delta(X)SSA(X)\\
&=0.
\end{align*}
\end{proof}
\end{subsection}
\end{section}

\begin{section}{%
Metric and Affine Conformal Geometry}\label{section:metricgeom}

We consider the metric extrinsic geometry of $f:M\to\RR^4$
in relation to the quantities associated to
\begin{align*}
L:=\aff{f\\1}:M\to\HH P^1.
\end{align*}
For brevity we write $<.,.>$ instead of $<.,.>_{\RR}$.

\begin{subsection}{Surfaces in Euclidean Space}\label{subsection:euclidean}

Let $N,R$ denote the left and right
normal vector of $f:M\to\HH$, i.e.
\begin{align*}
*df=Ndf=-dfR.
\end{align*}

\begin{proposition}\label{proposition:SFF} The
second fundamental form
$II(X,Y)=(X\cdot df(Y))^\perp$ of $f$ is given by
\begin{align}
II(X,Y)=\frac{1}{2}(*df(Y)dR(X)-dN(X)*df(Y)).\label{eq:SFF}
\end{align}
\end{proposition}
\begin{proof}
We know from Lemma~\ref{lemma:cstructure} that $v\mapsto N(x)vR(x)$ is an
 involution with the tangent space as its fixed point set:
\begin{align}
Ndf(Y)R=df(Y)\label{eq:tangfix}
\end{align}
Its $(-1)$-eigenspace is  the normal space, so we need to compute
\begin{align*}
II(X,Y)=\frac{1}{2}(X\cdot df(Y)-NX\cdot df(Y) R).
\end{align*}
But differentiation of \eqref{eq:tangfix} yields
\begin{align*}
dN(X)df(Y)R+NX\cdot df(Y) R+Ndf(Y)dR(X)=X\cdot df(Y),
\end{align*}
or
\begin{align*}
X\cdot df(Y)-NX\cdot df(Y) R
&=dN(X)df(Y)R+Ndf(Y)dR(X)\\
&=-dN(X)*df(Y)+*df(Y)dR(X).
\end{align*}
\end{proof}

\begin{proposition}
The mean curvature vector $\cH=\frac{1}{2}\trace II$ is given by
\begin{align}
\bar{\cH} df=\frac{1}{2}(*dR+RdR),\quad df\bar{\cH}=-\frac{1}{2}(*dN+NdN).
\label{eq:meanzwei}
\end{align}
\end{proposition}
\begin{proof}
By definition of the trace,
\begin{align}
4\cH |df|^2
&=*dfdR-dN*df-df*dR+*dNdf\\
&=-df(*dR+RdR)+(*dN+NdN)df,\label{eq:rllr}
\end{align}
but
\begin{align*}
(*dN+NdN)df&=*dNdf-dN*df=-dN\wedge df=-d(Ndf)\\
&=-df\wedge dR=-df(*dR+RdR).
\end{align*}
If follows that
\begin{align*}
2\cH |df|^2=-df(*dR+RdR),
\end{align*}
and
\begin{align*}
2\overline{\cH} df\overline{df}
=-(-*dR+dRR)\overline{df}
=(*dR+RdR)\overline{df}.
\end{align*}
Similarly for $N$.
\end{proof}

\begin{proposition} Let $K$ denote the Gaussian curvature  of
$(M,f^*<.,.>_{\RR})$
and let $K^\perp$ denote the normal curvature of $f$ defined by
\begin{align*}
K^\perp :=<R^\perp(X,JX)\xi,N\xi>_{\RR},
\end{align*}
where $X\in T_pM$, and $\xi\in \perp_p M$ are unit vectors. Then
\begin{align}
K|df|^2&=\frac{1}{2}(<*dR,RdR>+<*dN,NdN>)\label{eq:gauss}\\
K^\perp|df|^2&=\frac{1}{2}(<*dR,RdR>-<*dN,NdN>)\label{eq:nk}
\end{align}
\end{proposition}
\begin{proof}
\begin{align*}
K|df|^4(X)&=<II(X,X),II(JX,JX)>-|II(X,JX)|^2.
\end{align*}
Therefore

\begin{align*}
4K|df|^4
=&<*dfdR-dN*df,-df*dR+*dNdf>\\
&-<*df*dR-*dN*df,-dfdR+dNdf>\\
=&<N(dfdR+dNdf),-df*dR+*dNdf>\\
&-<N(df*dR+*dNdf),-dfdR+dNdf>\\
=&-<dfdR+dNdf,N(-df*dR+*dNdf)>\\
&<df*dR+*dNdf,N(-dfdR+dNdf)>\\
=&-<dfdR+dNdf,dfR*dR+N*dNdf>\\
&+<df*dR+*dNdf,dfRdR+NdNdf>\\
=&-<dfdR,dfR*dR>-<dfdR,N*dNdf>\\
&-<dNdf,dfR*dR>-<dNdf,N*dNdf>\\
&+<df*dR,dfRdR>+<df*dR,NdNdf>\\
&+<*dNdf,dfRdR>+<*dNdf,NdNdf>\\
=&-|df|^2<dR,R*dR>-<dfdR,N*dNdf>\\
&+<dNdf,Ndf*dR>-|df|^2<dN,N*dN>\\
&+|df|^2<*dR,RdR>+<df*dR,NdNdf>\\
&-<*dNdf,NdfdR>+|df|^2<*dN,NdN>\\
=&-|df|^2(<dR,R*dR>+<dN,N*dN>\\
&-<*dR,RdR>-<*dN,NdN>)\\
=&-2|df|^2(<dR,R*dR>+<dN,N*dN>).
\end{align*}
This proves the formula for $K$. Using \eqref{eq:SFF} and the Ricci equation
\begin{align*}
K^\perp=<N\,II(X,JX),II(X,X)-II(JX,JX)>,
\end{align*}
we find, after a similar computation,
\begin{align*}
4K^\perp|df|^2
&=<*dR-RdR,RdR>-<*dN-NdN,NdN>\\
&\quad +<df(*dR-RdR),NdNdf> -<(*dN-NdN)df,dfRdR>.
\end{align*}
On this we use \eqref{eq:rllr} to obtain \eqref{eq:nk}.
\end{proof}

As a corollary we have

\begin{proposition} The pull-back of the 2-sphere area under $R$ is given by
\begin{align*}
R^*dA=<*dR,RdR>.
\end{align*}
Integrating this for compact $M$ yields
\begin{align*}
\frac{1}{4\pi}\int_MK|df|^2=\frac{1}{2}(\deg R+\deg N).
\end{align*}
In 3-space ($R=N$) this is a version of the Gauss-Bonnet theorem.
\end{proposition}

\begin{proposition}\label{prop:wintaff} We obtain
\begin{align*}
(|\cH|^2-K-K^\perp)|df|^2=\frac{1}{4}|*dR-RdR|^2
\end{align*}
In particular, if $f:M\to\Im\HH=\RR^3$ then $K^\perp=0$, and the classical
Willmore integrand is given by
\begin{align}
(|\cH|^2-K)|df|^2=\frac{1}{4}|*dR-RdR|^2.\label{eq:willmoredrei}
\end{align}
\end{proposition}
\begin{proof}
Equations \eqref{eq:meanzwei}, \eqref{eq:gauss}, \eqref{eq:nk} give
\begin{align*}
(|\cH|^2-K-K^\perp)|df|^2
&=\frac{1}{4}|*dR+RdR|^2-<*dR,RdR>\\
&=\frac{1}{4}|*dR|^2+\frac{1}{4}|RdR|^2-\frac{1}{2}<*dR,RdR>\\
&=\frac{1}{4}|*dR-RdR|^2.
\end{align*}
\end{proof}

\end{subsection}

\begin{subsection}{The Mean Curvature Sphere in Affine
Coordinates}\label{subsection:mcsaffine}

We now discuss the characteristic properties of $S$ in affine coordinates.
We describe $S$ relative to the frame $\vect{1\\0},\vect{f\\1}$, i.e. we write
$S=GMG^{-1}$, where
\begin{align*}
G=\begin{pmatrix}1&f\\0&1\end{pmatrix}.
\end{align*}

First, $SL\subset L$ is equivalent to
$S:\HH^2\to\HH^2$ having the following matrix representation:
\begin{align}
S=
\begin{pmatrix}1&f\\0&1\end{pmatrix}
\begin{pmatrix}N&0\\-H&-R\end{pmatrix}
\begin{pmatrix}1&-f\\0&1\end{pmatrix}
\label{eq:Smatrix}
\end{align}
where $N,R,H:M\to\HH$. From $S^2=-I$
\begin{align}
N^2=-1=R^2, \quad RH=HN.\label{eq:NRH}
\end{align}
The choice of symbols is deliberate: $N$ and $R$ turn out to be
the left and right normal vectors of $f$, while $H$ is closely related to its
mean curvature vector
$\cH$.

The bundle $L$ has the nowhere vanishing section $\vect{f\\1}\in\Gamma(L)$.
Using
this section, we compute
\begin{align*}
*\delta\vect{f\\1}&=\pi\vect{*df\\0},\\
\delta S\vect{f\\1}
&=\pi d(S\vect{f\\1})=\pi
d(\vect{f\\1}(-R))=\pi(\vect{-dfR\\0}+\vect{f\\1}(-dR))=\pi\vect{-dfR\\0},\\
S\delta\vect{f\\1}&=\pi Sd\vect{f\\1}=\pi (\vect{Ndf\\0}+\vect{f\\1}(-Hdf))
=\pi \vect{Ndf\\0}.
\end{align*}
Therefore $*\delta=S\delta=\delta S$ is equivalent to
\begin{align*}
*df=Ndf=-dfR,
\end{align*}
and we have identified $N$ and $R$.

For the computation of the Hopf fields, we need $dS$. This is a
straight-forward but lengthy computation, somewhat simplified by the fact that
$GdG=dG=G^{-1}dG$. We skip the details and give the result:
\begin{gather*}
dS=G
\begin{pmatrix}
-dfH+dN&-dfR-Ndf\\
-dH&-dR+Hdf
\end{pmatrix}
G^{-1},\\
SdS=G
\begin{pmatrix}
-NdfH+NdN&0\\
HdfH+RdH-HdN&HdfR+RdR
\end{pmatrix}
G^{-1}.
\end{gather*}

From this we obtain
\begin{align*}
4Q&=SdS-*dS\\
&=G
\begin{pmatrix}
NdN-*dN&0\\
*dH+HdfH+RdH-HdN&2HdfR+RdR+*dR
\end{pmatrix}
G^{-1}
\\
4A&=SdS+*dS\\
&=G
\begin{pmatrix}
NdN+*dN-2NdfH&0\\
-*dH+HdfH+RdH-HdN&RdR-*dR
\end{pmatrix}
G^{-1}.
\end{align*}
The condition $Q|_L=0$, and the corresponding $AH\subset L$, which we have not
used so far, have the following equivalents:
\begin{align}
2Hdf=dR-R*dR,\label{eq:zweiHdf}\\
2dfH=dN-N*dN.\label{eq:zweidfH}
\end{align}

Together with equations \eqref{eq:meanzwei} we find
\begin{align*}
2Hdf&=dR-R*dR=-R(*dR+RdR)=-2R\bar{\cH} df,\\
2dfH&=dN-N*dN=-N(*dN+NdN)=2Ndf\bar{\cH}=-2dfR\bar{\cH},
\end{align*}
and therefore
\begin{align}
H=-\bar{\cH}N=-R\bar{\cH}.
\end{align}

\begin{remark}\label{remark:meancurvature}
Given an immersed holomorphic curve $L=\aff{f\\1}$, the
mean curvature vector of $f$ at $x\in M$ is determined by $S_x$. On the
other hand, $S_x$ is the mean curvature sphere of $S_x$,
see Example~\ref{example:tangSS}. Therefore $S_x$ and $f$ have, in fact, the
same mean curvature vector at $x$, justifying the name {\em mean curvature
sphere}.

\end{remark}

Equations \eqref{eq:zweiHdf}, \eqref{eq:zweidfH} simplify the coordinate
expressions for the Hopf fields, which we now write as follows

\begin{proposition}\label{proposition:w}
\begin{align}
4*Q&=G
\begin{pmatrix}
dN+N*dN&0\\
-2dH+w&0
\end{pmatrix}G^{-1},\label{eq:Qaffin}\\
4*A&=
G\begin{pmatrix}
0&0\\
w&dR+R*dR
\end{pmatrix}G^{-1},\label{eq:Aaffin}
\end{align}
where $G=\begin{pmatrix}1&f\\0&1\end{pmatrix}$, and
$w=dH+H*dfH+R*dH-H*dN$.

Using \eqref{eq:zweidfH} we can rewrite
\begin{align*}
w&=dH+R*dH+\frac{1}{2}H(NdN-*dN).
\end{align*}
\end{proposition}
\begin{proof} We only have to consider the reformulation of $w$. But
\begin{align*}
H*dfH-H*dN
=&\frac{1}{2}H*(dN-N*dN)-H*dN\\
=&-\frac{1}{2}H*(dN+N*dN)=\frac{1}{2}H(NdN-*dN).
\end{align*}
\end{proof}

\end{subsection}

\begin{subsection}{The Willmore Condition in Affine Coordinates}
\label{subsection:willmoreaffine}
We use the notations of the previous
Proposition~\ref{proposition:w}, and in addition abbreviate
\begin{align*}
v=dR+R*dR.
\end{align*}
Note that
\begin{align*}
\bar{v}=-dR+*dRR=-dR-R*dR=-v.
\end{align*}

\begin{proposition}\label{proposition:affinwillmoreint}
 The Willmore integrand is given by
\begin{align*}
<A\wedge*A>&=\frac{1}{16}|RdR-*dR|^2=\frac{1}{4}(|\cH|^2-K-K^\perp)|df|^2.
\end{align*}
For $f:M\to\RR^3$, this is the classical integrand
\begin{align*}
<A\wedge*A>=\frac{1}{4}(|\cH|^2-K)|df|^2.
\end{align*}
\end{proposition}
\begin{proof}
\begin{align*}
<A\wedge*A>&=\frac{1}{8}\trace_{\RR}(-A^2-(*A)^2)
=-\frac{1}{4}\trace_{\RR}(A^2)\\
&=-\frac{1}{4}\,4\Re (\frac{1}{4}v)^2= \frac{1}{16}|v|^2=
\frac{1}{16}|dR+R*dR|^2=
\frac{1}{16}|RdR-*dR|^2.
\end{align*}
Now see Proposition \ref{prop:wintaff} and, for the second
equality, \eqref{eq:willmoredrei}.
\end{proof}

We now express the Euler-Lagrange equation $d*A=0$ for Willmore surfaces in
affine coordinates. If we write $4*A=GMG^{-1}$, then
\begin{align*}
4d*A&=G(G^{-1}dG\wedge M+dM+M\wedge G^{-1}dG) G^{-1},
\end{align*}
and again using  $G^{-1}dG=dG$ we  easily find
\begin{align*}
4d*A=
G
\begin{pmatrix}
df\wedge w&df\wedge v\\
dw&dv+w\wedge df
\end{pmatrix}
G^{-1}.
\end{align*}
Most entries of this matrix vanish:

\begin{proposition}\label{proposition:vanishing}
We have
\begin{gather}
df\wedge w=0\label{eq:dfwedgew}\\
df\wedge v=0\label{eq:dfwedgev}\\
dv+w\wedge df=-(2dH-w)\wedge df=0.\label{eq:xxwedgedf}
\end{gather}
\end{proposition}

\begin{proof}
We have
\begin{align*}
df\wedge w
&=df\wedge dH +df\wedge R*dH +\frac{1}{2}df\wedge H(NdN-*dN)\\
&=df\wedge dH +df R\wedge *dH +\frac{1}{2}dfH\wedge (NdN-*dN)\\
&=\underbrace{df\wedge dH -*df \wedge *dH}_{=0}
 +\frac{1}{2}dfH\wedge (NdN-*dN),
\end{align*}
but
\begin{align*}
*df H&=df(-R)H=-df HN\\
*(NdN-*dN)&=(N*dN-N^2dN)=-N(NdN-*dN).
\end{align*}
Hence, by type, the second term vanishes as well, and we get
\eqref{eq:dfwedgew}.

A similar, but simpler, computation shows \eqref{eq:dfwedgev}

Next, using \eqref{eq:zweiHdf}, we consider
\begin{align*}
dv+w\wedge df&=d(dR+R*dR)+w\wedge df\\
&=d(-2Hdf)+w\wedge df\\
&=(-2dH+w)\wedge df\\
&=(\underbrace{-dH+R*dH}_{=:\alpha}
+\underbrace{\frac{1}{2}H(NdN-*dN)}_{\beta})\wedge df.
\end{align*}
Again we show $*\alpha=\alpha N,*\beta=\beta N$. Then
\eqref{eq:xxwedgedf} will follow by type.

Clearly
\begin{align*}
*(NdN-*dN)=N*dN+NdNN=(NdN-*dN)N,
\end{align*}
showing $*\beta=\beta N$.
Further
\begin{align*}
*\alpha-\alpha N
&=-*dH-RdH+dHN-R(*dH)N\\
&=-*dH-d(RH)+(dR)H+d(\underbrace{HN}_{=RH})-HdN-R*(d(\underbrace{HN}_{=RH})-HdN)
\\
&=+R^2*dH+(dR)H-HdN-R*((dR)H+RdH-HdN)\\
&=(dR)H-HdN-R*(dR)H+RH*dN)\\
&=(dR-R*dR)H-H(dN-N*dN)\\
&=2HdfH-H(2dfH)\\
&=0.
\end{align*}
\end{proof}

As a corollary we get:

\begin{proposition}\label{proposition:willoreaffine}
\begin{align*}
d*A =\frac{1}{4}G\begin{pmatrix}
0&0\\
dw&0
\end{pmatrix}G^{-1}=
\begin{pmatrix}
-fdw&-fdwf\\
dw&dwf
\end{pmatrix}.
\end{align*}
with
$w=dH+R*dH+\frac{1}{2}H(NdN-*dN)$.

Therefore $f$ is Willmore if and only if $dw=0$.
\end{proposition}

\begin{example}[Willmore Cylinder]
Let $\gamma:\RR\to\Im\HH$ be a unit-speed curve, and $f:\RR^2\to\HH$ the
cylinder defined by
\begin{align*}
f(s,t)=\gamma(s)+t
\end{align*}
with the conformal structure $J\frac{\partial}{\partial
s}=\frac{\partial}{\partial t}$. Then using  Proposition
\ref{proposition:willoreaffine}, we obtain, after some computation, that $f$ is
(non-compact) Willmore, if and only if
\begin{align*}
\frac{1}{2}\kappa^3+\kappa''-\kappa\tau^2=0,\quad (\kappa^2\tau)'=0.
\end{align*}
\end{example}
This is exactly the condition that $\gamma$ be a free elastic curve.
\end{subsection}
\end{section}

\begin{section}{%
Twistor Projections}

\begin{subsection}{Twistor Projections} \label{subsection:twistor}
Let $E\subset
H:=M\times\HH^2=M\times\CC^4$ be  a {\em complex} (not a quaternionic) line
subbundle over a Riemann surface $M$ with complex structure
$J_E$ induced from right multiplication by $i$ on $\HH^2$.

We define
$\delta_E\in\Omega^1(\Hom(E, H/E))$ by
\begin{align*}
\delta_E\phi&:=\pi_E d\phi,\quad\phi\in\Gamma(E),
\end{align*}
where $\pi_E: H\to H/E$ is the projection.

\begin{definition}
$E$ is called a {\em holomorphic curve} in $\CC P^3$, if
\begin{align*}
*\delta_E=\delta_E J_E.
\end{align*}
\end{definition}
This is equivalent to the fact that the holomorphic structure
\begin{align}
d''\psi=\frac{1}{2}(d\psi+i*d\psi)\label{eq:dzweistrich}
\end{align}
of $H$ maps $\Gamma(E)$ into itself, and hence induces a holomorphic
structure on the  complex line bundle  $E$.

A complex line bundle $E\subset H$ induces a quaternionic line bundle
\begin{align*}
L=E\HH=E\oplus Ej\subset H.
\end{align*}
The complex structure $J_E$ admits a unique extension to the structure of a
complex quaternionic bundle $(L,J)$, namely right-multiplication by $(-i)$
on $Ej$. Conversely, a complex quaternionic  line bundle $(L,J)\subset H$
induces a complex line bundle
\begin{align*}
E:=\{\phi\in L\,|\,J\phi=\phi i\}.
\end{align*}
\begin{definition}
We call $(L,J)$ the {\it twistor projection} of $E$, and
$E$ the {\it twistor lift} of
$(L,J)$.
\end{definition}
\begin{remark}
As in the quaternionic case, any map $f:M\to\CC P^3$ induces a complex line
bundle
$E$, where the fibre over $p$ is $f(p)$, and vice versa. Holomorphic
curves as defined above correspond to holomorphic curves in the sense of
complex analysis. The correspondence between $E$ and $(L,J)$ is mediated
by the Penrose twistor projection
$\CC P^3\to\HH P^1$.
\end{remark}

\begin{theorem}\label{theorem:twistorproj} Let $E\subset  H$ be a a complex
line subbundle over a Riemann surface $M$, and $(L,J)$ its twistor
projection.
\begin{enumerate}
\item Then (L,J) is a holomorphic curve, i.e.
\begin{align}
*\delta_L=\delta_L J,\label{eq:holcu}
\end{align}
if and only if
\begin{align*}
\frac{1}{2}(\delta_E+*\delta_E J_E)\in\Omega^1(\Hom(E,L/E))
\subset \Omega^1(\Hom(E,H/E)).
\end{align*}
In this case we have a differential operator
\begin{align*}
\tilde{D}:\Gamma(L)\to\Omega^1(L),
\psi\mapsto\tilde{D}\psi:=\frac{1}{2}(d\psi+*d(J\psi))
\end{align*}
Its $(1,0)$-part is given by
\begin{align}\label{eq:al}
A_L:=\frac{1}{2}(\tilde{D}+J\tilde{D}J)\in\Gamma(K\End_-(L)).
\end{align}
\item If $(L,J)$ is a holomorphic curve then
\begin{align*}
\frac{1}{2}(\delta_E+*\delta_E J_E)=\pi_EA_L|_E.
\end{align*}
Moreover,
\begin{align*}
\frac{1}{2}(\delta_E+*\delta_E J_E)=0\iff A_L=0.
\end{align*}
In other words: The twistor projections of holomorphic curves in
 $\CC P^3$ are exactly the
holomorphic curves in $\HH P^1$ with $A_L=0$.
\item Let $L$ be an immersed holomorphic curve with mean curvature sphere
congruence $S\in\Gamma(\End_-( H))$, and $J=S|_L$. Then
\begin{align*}
A=\frac{1}{4}(SdS+*dS)\in\Gamma(\bar{K}\End_-( H))
\end{align*}
satisfies
\begin{align*}
A|_L=A_L.
\end{align*}
\end{enumerate}
\end{theorem}
\begin{proof}
(i). If  $(L,J)$ is a holomorphic curve then, for any
$\psi\in\Gamma(L)$,
\begin{align*}
\frac{1}{2}\pi_L(d\psi+*d(J\psi))=0.
\end{align*}
But then
\begin{align*}
\frac{1}{2}(d\psi+*d(J\psi))\in \Omega^1(L)
\end{align*}
a fortiori for all $\psi=\phi\in\Gamma(E)$. It follows
\begin{align*}
\frac{1}{2}\pi_E(d\phi+*d (J_E\phi))\in
\Omega^1(L/E).
\end{align*}
Conversely,  $\frac{1}{2}\pi_E(d\phi+*d (J_E\phi))\in\Omega^1(L/E)$ for
$\phi\in\Gamma(E)$ implies
\begin{align*}
\frac{1}{2}(d\phi+*d(J_E\phi))\in \Omega^1(L),
\end{align*}
and therefore
\begin{align*}
*\delta_L|_E=\delta_L J|_E.
\end{align*}
Again for $\phi\in\Gamma(E)$
\begin{align*}
\frac{1}{2}(d(\phi j)+*d (J\phi j))
=\frac{1}{2}((d\phi) j+*d (J\phi) j))
=\frac{1}{2}(\underbrace{d\phi +*d (J_E\phi)}_{\in\Omega^1(L)})j
\in\Omega^1(L).
\end{align*}
This shows
\begin{align*}
*\delta_L=\delta_L J.
\end{align*}

By the preceding, $\tilde{D}$ maps into $\Omega^1(L)$. Its $(1,0)$-part is
\begin{align*}
\frac{1}{2}(\tilde{D}-J*\tilde{D}),
\end{align*}
but for $\psi\in\Gamma(L)$
\begin{align*}
*\tilde{D}\psi
=\frac{1}{2}(*d\psi-d(J\psi))
=-\tilde{D}J\psi.
\end{align*}
This proves \eqref{eq:al}.

(ii). For $\psi\in\Gamma(L)$ we have
\begin{align}
A_L\psi
&=\frac{1}{4}(d\psi+*d(J\psi)
+J(dJ\psi-*d\psi))\label{eq:cc}
\end{align}
But for $\phi\in\Gamma(E)$ we have $J(dJ\phi-*d\phi)
=J(d\phi +*d\phi i)i$, and hence
\begin{align*}
A_L\phi
&=\frac{1}{4}((d\phi+*d(J\phi))+J(d\phi +*d(J\phi))i).
\end{align*}
By assumption  $\frac{1}{2}(d\phi+*d(J\phi))$ has values in $L=E\oplus
Ej$, and $A_L\phi$ is its $Ej$-component, namely the component in
the  $(-i)$-eigenspace of $J|_L$. In particular,
\begin{align*}
\pi_EA_L\phi=\pi_E\frac{1}{2}(d\phi+*d(J\phi))
=\frac{1}{2}(\delta_E+*\delta_E J)\phi,
\end{align*}
and $\pi_E(A_L\phi)=0$ if and only
if $A_L\phi=0$. Since $A_L|_E$ determines $A_L$ by linearity,
$\frac{1}{2}\pi_E(d\psi+*d (\psi i))=0\iff A_L=0$.

(iii). For $\psi\in\Gamma(L)$
\begin{align*}
A\psi
&=\frac{1}{4}(SdS+*dS)\psi\\
&=\frac{1}{4}(S(d(S\psi)-Sd\psi)+*d(S\psi)-*Sd\psi)\\
&=\frac{1}{4}(S(d(S\psi)-*d\psi)+*d(S\psi)+d\psi).
\end{align*}
Comparison with \eqref{eq:cc} shows $A|_L=A_L$.
\end{proof}

\end{subsection}

\begin{subsection}{Super-Conformal Immersions.}

Given a surface conformally immersed into $\RR^4$, the image of a
tangential circle under the quadratic second fundamental form is  (a
double cover of) an ellipse in the normal space, centered at the mean
curvature vector, the so-called {\em curvature ellipse}. The surface is
called {\em super-conformal} if this ellipse is a circle.

If $N$ and $R$ are the left and right normal vector of $f$, then
according to Proposition~\ref{proposition:SFF}  we have
\begin{align*}
II(X,Y)=\frac{1}{2}(
*df(Y)dR(X)-dN(X)*df(Y)),
\end{align*}
and therefore
\begin{align*}
II(\cos\theta& X+\sin\theta JX,\cos\theta X+\sin\theta JX)
\\
=&\frac{1}{2}(*df(\cos\theta X+\sin\theta JX)dR(\cos\theta X+\sin\theta JX)
\\
 &-dN(\cos\theta X+\sin\theta JX)*df(\cos\theta X+\sin\theta JX))\\
=&\frac{1}{2}(df(\cos\theta JX-\sin\theta X)dR(\cos\theta X+\sin\theta JX)\\
 &-dN(\cos\theta X+\sin\theta JX)df(\cos\theta JX-\sin\theta X))\\
=&\frac{1}{2}(\cos^2\theta (df(JX)dR(X)-dN(X)df(JX))\\
 &\quad -\sin^2\theta(df(X)dR(JX)-dN(JX)df(X))\\
 &+\cos\theta\sin\theta(df(JX)dR(JX)-df(X)dR(X)\\
 &\qquad \qquad +dN(X)df(X)-dN(JX)df(JX)).
\end{align*}
Using $\cos^2\theta=\frac{1}{2}(1+\cos2\theta),
\sin^2\theta=\frac{1}{2}(1-\cos2\theta)$ we get
\begin{align*}
II(\cos\theta& X+\sin\theta JX,\cos\theta X+\sin\theta JX)
\\
=&\frac{1}{4}(
\underbrace{*df(X)dR(X)-dN(X)*df(X)}_{=2II(X,X)}
+\underbrace{*df(JX)dR(JX)-dN(JX)*df(JX)}_{=2II(JX,JX)})\\
 &+\frac{1}{4}\cos2\theta
 (df(JX)dR(X)-dN(X)df(JX)+df(X)dR(JX)-dN(JX)df(X))\\
&+\frac{1}{4}\sin2\theta(df(JX)dR(JX)-df(X)dR(X)+dN(X)df(X)-dN(JX)df(JX))\\
=&\cH|df(X)|^2\\
 &+\frac{1}{4}\cos2\theta
 (\underbrace{df(X)(*dR(X)-RdR(X))}_{=:a}
 -\underbrace{(*dN(X)-NdN(X))df(X)}_{=:b})\\
&+\frac{1}{4}\sin2\theta N(a+b).
\end{align*}

This is a circle if and only if $a-b$ and $N(a+b)$ are orthogonal and
have same length. This is clearly the case if $a=0$ or $b=0$, but
these are in fact the only possibilities. Assume that there exists
$P\in\HH, P^2=-1$ with
\begin{align}
N(a+b)=P(a-b),\label{eq:NP}
\end{align}
and note that
\begin{align*}
Na=aR,Nb=bR.
\end{align*}
We multiply \eqref{eq:NP} by $N$ from the left or by $R$ from the right to
obtain
\begin{align*}
-(a+b)=NP(a-b),\qquad -(a+b)=PN(a-b)
\end{align*}
respectively. Therefore $(PN-NP)(a-b)=0$, which implies $P=\pm N$, and
hence $a=0$ or $b=0$, or $a-b=0$. But then also $a+b=0$, whence $a=b=0$.

It follows that the immersion is super-conformal if and only if
\begin{align*}
*dR(X)-RdR(X)=0 \text{, or } *dN(X)-NdN(X)=0.
\end{align*}

By the preceding argument, this holds for a particular choice of $X$, but
 then it  obviously follows for all $X$.

We mention that $f\to\bar{f}$ exchanges $N$ and $R$, hence $f$
is super-conformal, if and only if $*dR-RdR=0$ for $f$ or for $\bar{f}$.
In view of
proposition~\ref{proposition:w}, this is equivalent to
$A|_L=0$, and by Theorem~\ref{theorem:twistorproj} we obtain:

\begin{theorem}\label{theorem:superconformal}
A conformally immersed Riemann surface $f:M\to\HH=\RR^4$
is super-conformal if and only if $\aff{f\\1}:M\to\HH P^1$ or
$\aff{\bar{f}\;\\1}:M\to\HH P^1$ is the twistor projection of a holomorphic
curve in $\CC P^3$.
\end{theorem}
\end{subsection}
\end{section}

\begin{section}{%
B\"acklund Transforms of Willmore Surfaces}

In this section we shall describe a method to construct new Willmore
surfaces from a given one. The construction depends on the choice of a
point $\infty$, and therefore generously offers a 4-parameter family of
such transformations. On the other hand, the necessary computations are not
invariant, and therefore ought to be done in affine coordinates.

The transformation theory is essentially local: This fact will be hidden
in the assumption that the transforms are again immersions. We shall also
ignore period problems.

\begin{subsection}{B\"acklund Transforms}
\label{subsection:onestep}
 Let
$f:M\to\HH$ be a Willmore surface with $N,R,H$, and
\begin{align*}
w=dH+H*dfH+R*dH-H*dN.
\end{align*}
Then
\begin{align*}
dw=0,
\end{align*}
and hence we can integrate it. Assume that $g:M\to \HH$ is an immersion
with
\begin{align}
dg=\frac{1}{2}w.\label{eq:defg}
\end{align}
{\small (Note that the integral of $w/2$ may have periods, so in general $g$
is defined only on a covering of $M$. We ignore this problem.)}

We want to show that $g$ is again a Willmore surface called a {\em
B\"acklund transform} of $f$.  Using this name, we refer to the fact that
in a given category of surfaces we construct new examples from old ones
by solving an ODE \eqref{eq:defg}, similar to the classical B\"acklund
transforms of K-surfaces, see Tenenblat~\cite{tenenblat}.

We denote the symbols associated to $g$ by a
subscript $(.)_g$, and want to prove $dw_g=0$.
The computation of $w_g$ can be done under the weaker assumption
\eqref{eq:dfwedgedg}, which holds in the case above, see
Proposition~\ref{proposition:vanishing}.

\begin{proposition}\label{proposition:willmoredaten}
Let $f,g:M\to\HH$ be immersions such that
\begin{align}
df\wedge dg=0.\label{eq:dfwedgedg}
\end{align}
Then $f$ and $g$ induce the same conformal structure on $M$, and
\begin{gather}
N_g=-R\label{eq:Ng},\\
dg(2dH_g-w_g)=-wdf.\label{eq:dgwg}
\end{gather}
\end{proposition}
\begin{proof}
Define $*$ using the conformal structure induced by $f$. Then
\begin{align*}
0=df\wedge dg=df*dg-df(-R)dg,
\end{align*}
which implies $*dg=-Rdg$. Hence $g$ is conformal, too, and $N_g=-R$.

For the next computations recall the
equations \eqref{eq:NRH}, and \eqref{eq:zweiHdf}, \eqref{eq:zweidfH}:
\begin{gather*}
HN=NR,\\
2dfH=dN-N*dN,\quad 2Hdf=dR-R*dR,\\
w=dH+H*dfH+R*dH-H*dN.
\end{gather*}
Then
\begin{align}
Rw
=&RdH+RH*dfH-*dH-RH*dN\nonumber\\
=&RdH+HN*dfH-*dH-HN*dN\nonumber\\
=&RdH-HdfH-*dH-H(N*dN-dN)-HdN\nonumber\\
=&RdH-HdN+HdfH-*dH.\label{eq:Rw1}
\end{align}
With $dRH+RdH=dHN+HdN$ this becomes
\begin{align}
Rw=dHN-*dH-dRH+HdfH.\label{eq:Rw2}
\end{align}
Next
\begin{align*}
2dgH_g=dN_g-N_g*dN_g=-dR-R*dR.
\end{align*}
Therefore
\begin{align*}
-dg\wedge dH_g
=\frac{1}{2}d(-dR-R*dR)=\frac{1}{2}d(dR-R*dR)=-dH\wedge df,
\end{align*}
or
\begin{align}
dg(*dH_g+R_gdH_g)=-(dHN-*dH)df.\label{eq:HgH}
\end{align}
We  now use \eqref{eq:Rw1} and \eqref{eq:HgH} to
compute
\begin{align*}
N_gdg&(2dH_g-w_g)\\
=&-dgR_g(2dH_g-w_g)\\
=&dg(-2R_gdH_g+R_gdH_g-H_gdN_g+H_gdgH_g-*dH_g)\\
=&-dg(R_gdH_g+*dH_g)+dgH_g(dgH_g-dN_g)\\
=&(dHN-*dH)df+dgH_g(dgH_g-dN_g)\\
=&(dHN-*dH)df+\frac{1}{4}(dN_g-N_g*dN_g)((dN_g-N_g*dN_g)-2dN_g)\\
=&(dHN-*dH)df-\frac{1}{4}(dR+R*dR)(dR-R*dR).
\end{align*}
Similarly, using \eqref{eq:Rw2},
\begin{align*}
-N_gwdf&=Rwdf\\
&=(dHN-*dH)df-(dR-Hdf)Hdf\\
&=(dHN-*dH)df-\frac{1}{4}(2dR-dR+R*dR)(dR-R*dR)\\
&=(dHN-*dH)df-\frac{1}{4}(dR+R*dR)(dR-R*dR).
\end{align*}
Comparison yields \eqref{eq:dgwg}.
\end{proof}

If $f$ is Willmore, and $g$ is defined by \eqref{eq:defg}, then
\begin{align*}
dg(2df+2dH_g-w_g)=2dgdf+dg(2dH_g-w_g)=(2dg-w)df=0.
\end{align*}
Hence
\begin{align}
w_g=2d(f+H_g),\label{eq:wg}
\end{align}
and $g$ is Willmore, too.

Now assume that $ h :=g-H$ is again an immersion. Then, by
Proposition~\ref{proposition:vanishing},
\begin{align*}
2dh \wedge df=(2dg-2dH)\wedge df=(w-2dH)\wedge df=0.
\end{align*}
Proposition~\ref{proposition:willmoredaten} applied to $(h,f)$  instead
of $(f,g)$ then says
\begin{align*}
-w_{h}dh &=df(2dH-w)
=df(2dH-2dg)
=-2dfdh .
\end{align*}
We find $w_{h}=2df$, whence $h$ is again a Willmore surface.
We call $g$ a {\em forward}, and $h$ a {\em backward B\"acklund
transform} of $f$. $h$ can be obtained without reference to $g$ by
integrating $d(g-H)=\frac{1}{2}w-dH$.

Note that $f$ is a forward B\"acklund transform of
$h$ because $df=\frac{1}{2}w_h$, and is also a backward transform of $g$
because  $df=\frac{1}{2}w_g-dH_g$, see \eqref{eq:wg}.

The concept of B\"acklund transformations depends on the
choice of affine coordinates. The following theorem
clarifies this situation.

\begin{theorem} Let
$L$ be a Willmore surface in
$\HH P^1$. Choose non-zero
$\beta\in(\HH^2)^*,a\in\HH^2$ such that $<\beta,a>=0$. Then
\begin{align*}
d<\beta,*Aa>=0=d<\beta,*Qa>.
\end{align*}
If $g,h:M\to\HH\subset\HH P^1$ are immersions that satisfy
\begin{align*}
dg=2<\beta,*Aa>,\quad dh=2<\beta,*Qa>,
\end{align*}
they are again Willmore surfaces, called {\em forward}
respectively {\em backward B\"acklund transforms} of $L$.
The free choice of $\beta$ implies that there is a whole $S^4$ of such
pairs of B\"acklund transforms. (Different choices of $a$ result in
Moebius transforms $g\to g\lambda$, or $h\to h\lambda$, for a constant
$\lambda$.)
\end{theorem}
\begin{proof}
Choose $b\in\HH^2,\alpha\in(\HH^2)^*$ such that $a,b$ and
 $\alpha,\beta$ are dual bases. Then
\begin{align*}
2<\beta,*Aa>=\frac{1}{2}w,\qquad 2<\beta,*Qa>=\frac{1}{2}w-dH,
\end{align*}
see Proposition~\ref{proposition:w}.
\end{proof}

We can now proceed from $g$ with another forward B\"acklund transform. To
do so, we must integrate $\frac{1}{2}w_g=d(f+H_g)$. But, up to a
translational constant, this yields
\begin{align}
\tilde{f}:=f+H_g.\label{eq:tf}
\end{align}
We now observe
\begin{lemma}
\begin{align*}
\vect{\tilde{f}\\1}\in\ker A.
\end{align*}
\end{lemma}
\begin{proof} Note that $\ker A=\ker *A$. By
Proposition~\ref{proposition:w} we have
\begin{align*}
4*A\vect{\tilde{f}\\1}
=&\begin{pmatrix}1&f\\0&1\end{pmatrix}
\begin{pmatrix}0&0\\w&dR+R*dR\end{pmatrix}
\begin{pmatrix}1&-f\\0&1\end{pmatrix}
\vect{f+H_g\\1}\\
=&\begin{pmatrix}1&f\\0&1\end{pmatrix}
\begin{pmatrix}0&0\\w&dR+R*dR\end{pmatrix}
\vect{H_g\\1}\\
=&\begin{pmatrix}1&f\\0&1\end{pmatrix}
\vect{0\\wH_g+\underbrace{dR+R*dR}_{=-dN_g+N_g*dN_g}}\\
=&\begin{pmatrix}1&f\\0&1\end{pmatrix}
\vect{0\\2dgH_g-2dgH_g}=0.
\end{align*}
\end{proof}
Similarly the twofold {\em backward} B\"acklund transform $\hat{f}$
satisfies
\begin{align*}
\vect{\hat{f}\\1}\HH\supset\im Q.
\end{align*}

But this means that away from the zeros of $A$ or $Q$ the 2-step
B\"acklund transforms of a Willmore surface $L$ in $\HH P^1$ can be
described simply as $\tL=\ker A$ or $\hL=\im Q$. In particular there are
no periods arising.

We obtain a chain of B\"acklund transforms
\begin{align*}
\begin{array}{ccccccccccccc}
\cdots&\to&\hat{f}&\to& h&\to&f&\to &g &\to &\tilde{f}& \to&\cdots\\
&&\parallel&&&&\parallel&& & &\parallel&&\\
\to&&\hL&&\to&&L&&\to&&\tL&&\to\\
\end{array}
\end{align*}
Of course, the chain may break down if we arrive at non-immersed
surfaces, or it may close up.
\end{subsection}

\begin{subsection}{Two-Step B\"acklund Transforms}

Let $L\subset H=M\times\HH^2$ be a Willmore surface, and assume
$A\not\equiv 0$ on each component of $M$. We want to describe directly the
two-step B\"acklund transform
$L\to\tL$, and compute its associated quantities (mean curvature sphere,
Hopf fields).

We state a fact about singularities that will be proved in the
appendix, see Section~\ref{subsection:supplement}.

\begin{proposition}
Let $L$ be a  Willmore surface in $\HH P^1$, and $A\not\equiv 0$ on each
component of $M$. Then there exists a unique line bundle $\tL\subset H$
such that on an open dense subset of $M$ we have:
\begin{align*}
\tL=\ker A,\text{ and } H=L\oplus\tL.
\end{align*}
A similar assertion holds for $\im Q$.
\end{proposition}

We shall assume that $\tL$ is immersed, and want to prove again that $\tL$ is
Willmore.

\begin{theorem}\label{theorem:tQistA}
For the 2-step B\"acklund transform $\tL$ of $L$ we have
\begin{align}
\tQ=A.\label{eq:tQistA}
\end{align}
Hence $\tL$ is again a Willmore surface.
\end{theorem}

Let $\tS,\tdelta,\tQ$, etc. denote the operators associated with $\tL$.
\begin{lemma}\label{lemma:tdelta}
\begin{align*}
*\tdelta=-S\tdelta.
\end{align*}
\end{lemma}
\begin{proof}
Since $A|_{\tL}=0$ we interpret $A\in\Omega^1(\Hom(H/\tL,H))$.
On a dense open subset of $M$ then $A(X):H/\tL\to H$ is injective for
any $X\neq0$.

For $\phi\in\Gamma(\tL)$ we get
\begin{align*}
0=&d(*A)\phi
=d(\underbrace{*A\phi}_{=0})+*A\wedge d\phi=*A*d\phi+A d\phi\\
=&-AS*\tdelta\phi+A\tdelta\phi
=-AS(*\tdelta+S\tdelta)\phi.
\end{align*}
The injectivity of $A$ then proves the lemma.
\end{proof}

\begin{proof}[Proof of the theorem]
Motivated by the lemma, we
relate $\tS$ to $-S$ rather than to
$S$. We put
\begin{align*}
\tS=:-S+B.
\end{align*}
Then
\begin{align*}
4\tQ&=\tS d\tS-*d\tS\\
&=Bd\tS-(Sd\tS+*d\tS)\\
&=Bd\tS-(SdB+*dB)+(SdS+*dS)\\
&=4A+Bd\tS-(SdB+*dB).
\end{align*}
The proof will be completed with the following lemma which shows that
$\tQ$ -- like $A$ -- has values in $L$, while the ``B-terms'' take values
in $\tL$.
\end{proof}

\begin{lemma} We have
\begin{gather}
\im B\subset\tL,\label{eq:imageB}\\
\im (*dB+SdB)\subset \tL,\label{eq:imagedB}\\
L\subset\ker B,\label{eq:kerB}\\
\im \tQ\subset L.\label{eq:imtQ}
\end{gather}
\end{lemma}

\begin{proof}
Recall that $\tL$ is $S$-stable. It
is of course also $\tS$-stable, and therefore
\begin{align}
B\tL\subset \tL. \label{eq:Bstable}
\end{align}
Now $\tL$ is immersive, and therefore
$\im\tdelta=H/\tL$. Thus \eqref{eq:imageB} will follow if we can show
$\tilde{\pi}Bd\phi=0$ for $\phi\in\Gamma(\tL)$. But, using
Lemma~\ref{lemma:tdelta},
\begin{align*}
\tilde{\pi}Bd\phi
&=\tilde{\pi}Sd\phi+\tilde{\pi}\tS d\phi
=S\tilde{\pi}d\phi+\tS\tilde{\pi}d\phi
=S\tdelta\phi+\tS\tdelta\phi\\
&=-*\tdelta\phi+*\tdelta\phi=0.
\end{align*}

Next, for $\chi\in\Gamma(H)$ we have
\begin{align*}
\tilde{\pi}(*dB+SdB)\chi
&=\tilde{\pi}(*d(B\chi)+Sd(B\chi)-\underbrace{B*d\chi
-SBd\chi}_{\tL-\text{valued}})\\ &=(*\tdelta +S\tdelta)B\chi\\
&=0.\qquad (\text{Lemma \ref{lemma:tdelta}})
\end{align*}
This proves \eqref{eq:imagedB}.

On the other hand, for $\psi\in\Gamma(L)$,
\begin{align*}
\tilde{\pi}(*dB-SdB)\psi
&=\tilde{\pi}\underbrace{(*dS-SdS)\psi}_{=-4Q\psi=0}
+\tilde{\pi}(*d\tS-Sd\tS)\psi\\
&=\tilde{\pi}\underbrace{(*d\tS+\tS
d\tS)\psi}_{=4 \tA \psi\in\Gamma(\tL)}
-\tilde{\pi}\underbrace{(Bd\tS)\psi}_{\in\Gamma(\tL)}\\
&=0.
\end{align*}
Together with the previous equation we obtain
 $\tilde{\pi}dB|_L=0$, and,  for  $\psi\in\Gamma(L)$,
\begin{align*}
\tdelta B\psi=\tilde{\pi}(d(B\psi))=\tilde{\pi}((dB)\psi-Bd\psi)=
\tilde{\pi}dB\psi=0.
\end{align*}
But $\tL$ is an immersion, and therefore $B\psi=0$, proving
\eqref{eq:kerB}.

Finally, for $\psi\in\Gamma(L)$,
\begin{align*}
4\tQ\psi
&=\tS d\tS\psi-*d\tS\psi\\
&=\tS d\tS\psi-d\psi+\tS*d\psi-(-d\psi+\tS*d\psi+*d\tS\psi)\\
&=\tS (d\tS\psi+\tS d\psi+*d\psi)-*(*d\psi+\tS d\psi+ d\tS\psi)\\
&=(\tS-*)(d(\tS\psi)+*d\psi)\\
&=-(\tS-*)(d(S\psi)-*d\psi)\quad\text{using }\eqref{eq:kerB}.
\end{align*}
But $\pi(d(S\psi)-*d\psi)=(\delta S-*\delta)\psi=0$. So
$d(S\psi)-*d\psi\in\Gamma(L)$, and this is stable under $\tS=B-S$.
Therefore $ \tQ L\subset L$. Since $ \tQ \tL=0$, this proves
\eqref{eq:imtQ}.
\end{proof}

Taking the two-step backward transform of $\tL$, we get $\im \tQ=\im A=L$.
Hence $\hat{\tL}=L$. We remark that the results of this section  similarly
apply to the backward two-step B\"acklund transformation
$L\to\hL=\im Q$. As a corollary of \eqref{eq:tQistA} and its analog
$\hat{A}=Q$ we obtain

\begin{theorem}
\begin{align*}
\hat{\tL}=L=\tilde{\hL}.
\end{align*}
\end{theorem}
\end{subsection}
\end{section}

\begin{section}{%
Willmore Surfaces in $S^3$}

Let $<.,.>$ be an indefinite hermitian inner product on $\HH^2$. To be
specific, we choose
\begin{align*}
<v,w>:=\bar{v_1}w_2+\bar{v_2}w_1.
\end{align*}
Then the set of isotropic lines $<l,l>=0$ defines an $S^3\subset \HH P^1$,
while  the complementary 4-discs are hyperbolic 4-spaces, see
Example~\ref{example:hyperbolic}. We have
\begin{align}
\begin{pmatrix}
a&b\\c&d
\end{pmatrix}^*
=
\begin{pmatrix}
\bar{d}&\bar{b}\\\bar{c}&\bar{a}
\end{pmatrix},
\label{eq:adjoint}
\end{align}
and the same holds for matrix representations with respect to a basis
$(v,w)$ such that
\begin{align*}
<v,v>=0=<w,w>,\quad <v,w>=1.
\end{align*}

\begin{subsection}{Surfaces in $S^3$.}
Let $L$ be an isotropic line bundle with mean curvature sphere $S$. We
look at the adjoint map $M\to\cZ,p\mapsto S^*_p$ with respect to
 $<.,.>$. Clearly $S^*$ stabilizes $L^\perp$, and
$L=L^\perp$ implies
\begin{align*}
S^*L=S^*L^\perp=L^\perp=L.
\end{align*}
Similarly,
\begin{align*}
(dS^*)L=(dS)^*L^\perp\subset L^\perp=L.
\end{align*}
Moreover, if $Q^\dagger$ belongs to $S^*$, then
\begin{align*}
Q^\dagger&=\frac{1}{4}(S^*dS^*-*dS^*)\\
&=\frac{1}{4}(dS S-*dS)^*\\
&=-\frac{1}{4}(SdS+*dS)^*\\
&=-A^*.
\end{align*}
Therefore $\ker Q^\dagger=(\im (Q^\dagger)^*)^\perp=(\im A)^\perp\supset
L^\perp=L$.

By the uniqueness of the mean curvature sphere,
see Theorem~\ref{theorem:mcs}, it follows that
$S^*=S$. Conversely, if $S^*=S$  and
$S\psi=\psi\lambda$, then
\begin{align*}
\bar{\lambda}<\psi,\psi>=<S\psi,\psi>=<\psi,S\psi>
=<\psi,\psi>\lambda=\lambda<\psi,\psi>.
\end{align*}
Now $S^2=-I$ implies $\lambda^2=-1$, and therefore we get $<\psi,\psi>=0$.
\begin{proposition}
An immersed holomorphic curve $L$ in $\HH P^1$ is isotropic, i.e. a
surface in $S^3$, if and only if
$S=S^*$.
\end{proposition}

\end{subsection}

\begin{subsection}{Hyperbolic 2-Planes}
In the half-space or Poincar\'e model of the hyperbolic space, geodesics
are euclidean circles that orthogonally intersect the boundary. We
consider the models of hyperbolic 4-space in $\HH P^1$, and want to
identify their totally geodesic hyperbolic 2-planes, i.e. those 2-spheres
in $\HH P^1$ that orthogonally intersect the separating isotropic $S^3$.
Using the affine coordinates, from Example~\ref{example:hyperbolic}, we
consider the reflexion
$\HH\to\HH,x\mapsto -\bar{x}$ at $\Im\HH=S^3$. This preserves either of
the metrics given in the examples of Section
\ref{subsection:metrics}. In particular, it induces an isometry of the
standard  Riemannian metric of $\HH P^1$ which fixes
$S^3$. Given a 2-sphere $S\in\End(\HH^2), S^2=-I$, that intersects $S^3$ in
a point $l$, we use affine coordinates, as in
Example~\ref{example:hyperbolic}, with $l=v\HH$ and $w$ such
that
\begin{align*}
<v,v>=<w,w>=0,<v,w>=1.
\end{align*}
Then
\begin{align*}
S=
\begin{pmatrix}N&-H\\0&-R\end{pmatrix}
\end{align*}
with $N^2=R^2=-1,HN=RH$,
and $S'\subset\HH$ is the locus of
\begin{align*}
Nx+xR=H.
\end{align*}
If $S'$ is invariant under the reflexion at
$S^3$, then it also is the locus of $-N\bar{x}-\bar{x}R=H$ or
\begin{align*}
Rx+xN=\bar{H}.
\end{align*}
According to Section~\ref{subsection:twospheres},
the triple $(H,N,R)$ is unique up to sign. This implies either
\begin{align*}
(H,N,R)=(\bar{H},R,N) \text{ or } (H,N,R)=(-\bar{H},-R,-N).
\end{align*}
By \eqref{eq:adjoint} either $S^*=S$, and the 2-sphere lies within the
3-sphere, or it intersects orthogonally, and $S^*=-S$.  We summarize:

\begin{proposition} A 2-sphere $S\in\cZ$ intersects the hyperbolic
4-spaces determined by an indefinite inner product in hyperbolic
2-planes if and only if
$S^*=-S$.
\end{proposition}
\end{subsection}

\begin{subsection}{Willmore Surfaces in
$S^3$ and Minimal Surfaces in Hyperbolic 4-Space}

Let $L$ be a connected Willmore surface in $S^3\subset\HH P^1$, where $S^3$
is the isotropic set of an indefinite hermitian form on $\HH^2$.
Then its mean curvature sphere satisfies
\begin{align*}
S^*=S.
\end{align*}
Let us assume that $A\not\equiv0$, and let $\tL=\ker A$ and
$\hL=\im Q$ be the 2-step B\"acklund transforms of
$L$.
\begin{lemma}
\begin{align*}
\hL=\tL.
\end{align*}
\end{lemma}
\begin{proof}
First we have
\begin{align}
Q^*
&=\frac{1}{4}(SdS-*dS)^*
=\frac{1}{4}(dSS-*dS)\nonumber\\
&=\frac{1}{4}(-SdS-*dS)
=-A.\label{eq:qstern}
\end{align}
Now $\hL=\im Q$ is $S$-stable, and $S^*=S$ and $S\phi=\phi\lambda$ imply
$<\phi,\phi>=0$. Therefore $<\hL,\hL>=0$, and on  a dense open subset of
$M$
\begin{align*}
\hL=\hL^\perp=(\im Q)^\perp=\ker Q^*=\ker A=\tL.
\end{align*}
\end{proof}

\begin{lemma}\label{lemma:tS}
\begin{align*}
\tS=-S
\end{align*}
for the mean curvature sphere $\tS$ of $\tL$.
\end{lemma}
\begin{proof}
First $\tL=\hL$ is obviously $(-S)$-stable. It is trivially invariant under
$A$ and $Q$ and, therefore, under $d(-S)=2(*A-*Q)$. Finally, the $Q$ of
$(-S)$ is
\begin{align*}
\frac{1}{4}((-S)d(-S)-*d(-S))=A,
\end{align*}
and this vanishes on $\tL$.
The unique characterization of the mean curvature sphere by these three
properties implies $\tS=-S$.
\end{proof}
We now turn to the 1-step B\"acklund transform of $L$. If $dF=2*A$,
then
\begin{align*}
d(F+F^*)=2*A+2*A^*\underset{\eqref{eq:qstern}}{=}2*A-2*Q=-dS.
\end{align*}
Because $S^*=S$, we can choose suitable initial conditions for $F$ such
that
\begin{align}
F+F^*=-S.\label{eq:normierung}
\end{align}
We now use affine coordinates with $L=\aff{f\\1}$. Then the lower
left entry $g$  of $F$ is a B\"acklund transform of $f$, and
\eqref{eq:Smatrix} and
\eqref{eq:normierung} imply
\begin{align*}
g+\bar{g}=H.
\end{align*}
We want to compute the mean curvature sphere $S_g$. From the properties of
B\"acklund transforms we know
\begin{align}
N_g=-R,\quad H_g=\tilde{f}-f,
\end{align}
see \eqref{eq:Ng}, \eqref{eq:tf}. Likewise, $\tilde{N}=-R_g$. From
Lemma~\ref{lemma:tS} we obtain
\begin{align*}
\begin{pmatrix}1&f\\0&1\end{pmatrix}&
\begin{pmatrix}-N&0\\H&R\end{pmatrix}
\begin{pmatrix}1&-f\\0&1\end{pmatrix}
=
\begin{pmatrix}1&\tilde{f}\\0&1\end{pmatrix}
\begin{pmatrix}\tilde{N}&0\\-\tilde{H}&-\tilde{R}\end{pmatrix}
\begin{pmatrix}1&-\tilde{f}\\0&1\end{pmatrix}\\
&=
\begin{pmatrix}1&f\\0&1\end{pmatrix}
\begin{pmatrix}1&H_g\\0&1\end{pmatrix}
\begin{pmatrix}\tilde{N}&0\\-\tilde{H}&-\tilde{R}\end{pmatrix}
\begin{pmatrix}1&-H_g\\0&1\end{pmatrix}
\begin{pmatrix}1&-f\\0&1\end{pmatrix}\\
&=
\begin{pmatrix}1&f\\0&1\end{pmatrix}
\begin{pmatrix}\tilde{N}-H_g\tilde{H}&*\\ -\tilde{H}&*\end{pmatrix}
\begin{pmatrix}1&-f\\0&1\end{pmatrix}.
\end{align*}
This implies $H=-\tilde{H}$ and $-N=\tilde{N}-H_g\tilde{H}$, whence
\begin{align*}
-R_g=\tilde{N}=-N+(f-\tilde{f})H.
\end{align*}
In particular $f-\tilde{f}\in\Im\HH$, since $H=0$ on an open set would mean
$w=0$ on that set. It follows that
\begin{align*}
S_g=\begin{pmatrix}1&g\\0&1\end{pmatrix}
\begin{pmatrix}-R&0\\f-\tilde{f}&-N+(f-\tilde{f})H\end{pmatrix}
\begin{pmatrix}1&-g\\0&1\end{pmatrix},
\end{align*}
and, because $R=N$ and $H\in\RR$ for $f:M\to\Im\HH=\RR^3$,
\begin{align*}
S_g^*
&=\begin{pmatrix}1&g-H\\0&1\end{pmatrix}
\begin{pmatrix}\overline{-N+(f-\tilde{f})H}&0\\
\overline{f-\tilde{f}}&-\bar{R}\end{pmatrix}
\begin{pmatrix}1&H-g\\0&1\end{pmatrix}\\
&=\begin{pmatrix}1&g-H\\0&1\end{pmatrix}
\begin{pmatrix}N+(\tilde{f}-f)H&0\\
\tilde{f}-f&N\end{pmatrix}
\begin{pmatrix}1&H-g\\0&1\end{pmatrix}\\
&=\begin{pmatrix}1&g\\0&1\end{pmatrix}
\begin{pmatrix}N&0\\
\tilde{f}-f&N+(\tilde{f}-f)H\end{pmatrix}
\begin{pmatrix}1&-g\\0&1\end{pmatrix}\\
&=-S_g.
\end{align*}
We have now shown that the mean curvature spheres of $g$ intersect $S^3$
orthogonally, and  therefore are hyperbolic planes. We know that, using
affine coordinates and a Euclidean metric, the mean curvature spheres are
tangent to $g$ and have the same mean curvature vector as $g$. This
property remains under conformal changes of the ambient metric. Therefore,
in the hyperbolic metric, $g$ has mean curvature 0, and hence is
minimal. If $A\equiv0$, then $w=0$, and the ``B\"acklund transform'' is
constant, which may be considered as a degenerate minimal surface. In
general $g$ will be singular in the (isolated) zeros of $dg=\frac{1}{2}w$,
but minimal elsewhere.

We show the converse: Let $L$ be an immersed holomorphic curve,
minimal in hyperbolic 4-space, i.e. with $S^*=-S$. Then
\begin{align*}
A^*=\frac{1}{4}(SdS+*dS)^*
   =\frac{1}{4}(dSS-*dS)
   =-\frac{1}{4}(SdS+*dS)=-A,
\end{align*}
and therefore also
\begin{align*}
(d*A)^*=-d*A.
\end{align*}
From Proposition~\ref{proposition:willoreaffine} we have
\begin{align*}
d*A=
\begin{pmatrix}
-fdw&-fdwf\\
dw&dwf
\end{pmatrix}.
\end{align*}
Therefore
\begin{gather*}
dw=-\overline{dw},\qquad \overline{fdw}=dwf,
\end{gather*}
and hence
\begin{align*}
dw(f+\bar{f})=0.
\end{align*}
But $f$ is not in $S^3$, and therefore $dw$=0, i.e $L$ is Willmore.
Similarly,  Proposition~\ref{proposition:w} yields
\begin{align*}
*A=\begin{pmatrix}*&*\\w&*\end{pmatrix},
\end{align*}
and $A^*=-A$ implies $w=-\bar{w}$. From $S*=-S$ we know $\bar{H}=-H$,
and the backward B\"acklund transform $h$ with $dh=\frac{1}{2}-dH$ and
suitable initial conditions is in $\Im\HH=\RR^3$.

To summarize
\begin{theorem}[Richter \cite{richter}]
Let  $<.,.>$ be  an indefinite hermitian product on $\HH^2$. Then
the isotropic lines form an $S^3\subset\HH P^1$, while the two
complementary discs inherit complete hyperbolic metrics. Let $L$ be a
Willmore surface in $S^3\subset\HH P^1$. Then a suitable forward B\"acklund
transform of $L$ is hyperbolic mini\-mal. Conversely, an immersed
holomorphic curve that is  hyperbolic minimal is Willmore, and a
suitable backward B\"acklund transformation is a Willmore surface in
$S^3$. (In both cases the B\"acklund transforms may have singularities.)
\end{theorem}
\end{subsection}
\end{section}

\begin{section}{%
Spherical Willmore Surfaces in  $\HH P^1$}\label{section:montiel}

In this section we sketch a proof of the following theorem of Montiel,
which generalizes an earlier result of Bryant \cite{bryant} for Willmore
spheres in $S^3$.

\begin{theorem}[Montiel \cite{montiel}]\label{theorem:montiel}
A Willmore {\em sphere} in
$\HH P^1$ is a twistor projection of a holomorphic or anti-holomorphic
curve in $\CC P^3$, or, in suitable affine coordinates,
corresponds to a minimal surface in $\RR^4$.
\end{theorem}

The material differs from what we have treated so far:
The theorem is global, and therefore requires global methods of proof.
These are imported from complex function theory.

\begin{subsection}{Complex Line Bundles: Degree and Holomorphicity}
Let $E$ be a complex vector bundle. We keep the symbol
$J\in\End(H)$ for the endomorphism given by multiplication with the
imaginary unit $i$.

We denote by $\bar{E}$ the bundle where $J$ is replaced by $-J$. If $<.,.>$ is
a hermitian metric on $E$, then
\begin{align*}
\bar{E}\to E^*=E^{-1},\psi\to<\psi,.>
\end{align*}
is an isomorphism of complex vector bundles. Also note that for complex
{\em line} bundles $E_1,E_2$ the bundle $\Hom(E_1,E_2)$ is again a complex
line bundle.

There is a powerful integer invariant for complex line bundles $E$ over a
compact Riemann surface: the {\em degree}. It classifies these bundles up
to isomorphism. Here are two equivalent definitions for the degree.
\begin{itemize}
\item Choose a hermitian metric $<.,.>$ and a compatible connection
$\nabla$ on $E$. Then $<R(X,Y)\psi,\psi>=0$ for the curvature tensor $R$
of $\nabla$. Therefore $R(X,Y)=-\omega(X,Y)J$ with a real
2-form $\omega\in\Omega^2(M)$. Define
\begin{align*}
\deg(E):=\frac{1}{2\pi}\int_M\omega.
\end{align*}
\item Choose a section $\psi\in\Gamma(E)$ with isolated zeros. Then
\begin{align*}
\deg(E):=\ord\phi:=\sum_{\phi(p)=0}\ind_{p}\phi.
\end{align*}
The index of a zero $p$ of $\phi$ is defined using a local non-vanishing
section $\psi$ and a holomorphic parameter $z$ for $M$ with $z(0)=p$. Then
$\phi(z)=\psi(z)\lambda(z)$ for some complex function $\lambda:\CC\subset
U\to\CC$ with isolated zero at 0, and
\begin{align*}
\ind_{p}\phi=\frac{1}{2\pi i}\int_{\gamma}\frac{dz}{\lambda(z)},
\end{align*}
where $\gamma$ is a small circle around $0$.
\end{itemize}
We state fundamental properties of the degree.  We have
\begin{gather*}
\deg(\bar{E})=\deg{E^{-1}}=-\deg E,\\
\deg\Hom(E_1,E_2)=-\deg E_1+\deg E_2.
\end{gather*}
More generally,
\begin{align*}
\deg(E_1\otimes E_2)=\deg E_1+\deg E_2.
\end{align*}
\begin{example}
Let $M$ be a compact Riemann surface of genus $g$, and $E$ its tangent
bundle, viewed as a complex line bundle. We compute its degree using the
first definition. The curvature tensor of a surface with Riemannian metric
$<.,.>$ is given by
$R(X,Y)=K(<Y,.>X-<X,.>Y)$, where $K$ is the Gaussian curvature. Let $Z$ be
a (local) unit vector field and $<.,>$ compatible with $J$. Then
\begin{align*}
\omega(X,Y)&=\frac{1}{2}\trace_{\RR} R(X,Y)J\\
&=\frac{K}{2}(<Y,JZ><X,Z>-<X,JZ><Y,Z>\\
&\qquad-<Y,Z><X,JZ>+<X,Z><Y,JZ>)\\
&=K(<Y,JZ><X,Z>-<X,JZ><Y,Z>)\\
&=K\det\begin{pmatrix}<X,Z>&<X,JZ>\\<Y,Z>&<Y,JZ>)\end{pmatrix}\\
&=K\,dA(X,Y).
\end{align*}
We integrate this using Gauss-Bonnet, and find
$2\pi\chi(M)=2\pi(2-2g)=2\pi\deg(E)$.
For the canonical bundle
\begin{align*}
K:=E^{-1}=\Hom(TM,\CC)=\{\omega\in\Hom_\RR(TM,\CC)\,|\,\omega(JX)=i\omega(X)\}
\end{align*}
we therefore find
\begin{align*}
\deg(K)=2g-2.
\end{align*}
\end{example}

\begin{definition}
Let $E$ be a complex vector bundle. A {\em holomorphic structure} for $E$
is a complex linear map a map $\dbar$ from the sections of $E$ into
the $E$-valued complex anti-linear 1-forms $\bar{K}E$
\begin{align*}
\dbar:\Gamma(E)\to\Gamma(\bar{K}E)
\end{align*}
satisfying
\begin{align*}
\dbar(\lambda\psi)=(\dbar\psi)\lambda+\psi(\dbar\lambda).
\end{align*}
Here $\dbar\lambda:=\frac{1}{2}(d\lambda+i*d\lambda)$.
(Local) sections $\psi\in\Gamma(E|_U)$ are called {\em holomorphic}, if
$\dbar\psi=0$. We denote by $H^0(E|_U)$ the vector space of holomorphic
sections over $U$.
\end{definition}

If $E$ is a complex {\em line} bundle with holomorphic structure, and
$\psi\in H^0(E)\backslash\{0\}$, then the zeros of $\psi$ are isolated and
of positive index because holomorphic maps preserve orientation. In
particular, if $M$ is compact and $\deg E<0$, then any global
holomorphic section in $E$ vanishes identically.

In the proof of the Montiel theorem we shall apply the concepts of degree
and holomorphicity to several complex bundles obtained from quaternionic
ones. We relate these concepts.

\begin{definition}
If $(L,J)$ is a complex {\em quaternionic} line bundle, then
\begin{align*}
E_L:=\{\psi\in L\,|\,J\psi=\psi i\}
\end{align*}
is a complex line bundle. We define
\begin{align*}
\deg L:=\deg E_L.
\end{align*}
\end{definition}

\begin{lemma} If $L_1,L_2$ are complex quaternionic line bundles, and
$E_i:=E_{L_i}$, then
\begin{align*}
\Hom_+(L_1,L_2)&\to\Hom_{\CC}(E_1,E_2)\\
B&\mapsto B|_{E_1}
\end{align*}
is an isomorphism of complex vector bundles. In particular
\begin{align*}
\deg\Hom_+(L_1,L_2)=-\deg L_1+\deg L_2.
\end{align*}
\end{lemma}
The proof is straightforward. We now discuss one  example in detail.

\begin{example}\label{example:Aholo} We consider an immersed holomorphic
curve
\begin{align*}
L\subset H=M\times\HH^2
\end{align*}
in $\HH P^1$ with mean curvature sphere $S$. The bundle
$K\End_-(H)$ is a complex vector bundle, the complex structure being given
by post-composition with $S$. For
$B\in\Gamma(K\End_-(H))$ we define
\begin{align*}
(\dbar_XB)(Y)\psi=\dbar_X(B(Y)\psi)-B(\dbar_XY)\psi-B(Y)\partial_X\psi,
\end{align*}
where
\begin{gather*}
\dbar_XY:=\frac{1}{2}([X,Y]+J[JX,Y]),\\
\dbar\psi=\frac{1}{2}(d+S*d)\psi,\quad
\partial\psi=\frac{1}{2}(d-S*d)\psi\text{ for }\psi\in\Gamma(H).
\end{gather*}
Direct computation shows that this is in fact a holomorphic
structure, namely that induced on
\begin{align*}
K\End_-(H)=K\Hom_+(\bar{H},H)=K\Hom_{\CC}(\bar{H},H)
\end{align*}
by $\dbar$ on $TM$, and the above (quaternionic) holomorphic
structures $\dbar$ on $H$ and  $\partial$ on $\bar{H}$.

\begin{lemma}
\begin{align*}
(d*A)(X,JX)=-2(\dbar_X A)(X).
\end{align*}
\end{lemma}
\begin{proof}
Let $X$ be a local holomorphic vector field, i.e.  $[X,JX]=0$, see
Remark~\ref{remark:holoVF}, and $\psi\in\Gamma(H)$. Then
\begin{align*}
(d*A)(X,JX)\psi
&=(-X\cdot A(X)-(JX)\cdot SA(X)-A(\underbrace{[X,JX]}_{=0})\psi\\
&=-(d(\underbrace{A(X)\psi}_{=:\phi})+*d(SA(X)\psi))(X)\\
&\quad +A(X)d\psi(X)+SA(X)*d\psi(X)\\
&=-(d\phi+*d(S\phi))(X)
 +A(X)(d\psi-S*d\psi)(X).
\end{align*}
Now
\begin{align*}
d\phi+*d(S\phi)
&=(\partial+\dbar+A+Q)\phi+*(\partial+\dbar+A+Q)S\psi\\
&=(\partial+\dbar+A+Q)\phi+(S\partial-S\dbar+SA-SQ)S\psi\\
&=(\partial+\dbar+A+Q)\phi+(-\partial+\dbar+A-Q)\psi\\
&=2(\dbar+A)\phi\\
&=2\dbar(A(X)\psi)+2AA(X)\psi.
\end{align*}
Similarly
\begin{align*}
d\psi-S*d\psi
&=(\partial+\dbar+A+Q)\psi-S*(\partial+\dbar+A+Q)\psi\\
&=(\partial+\dbar+A+Q)\psi-S(S\partial-S\dbar+SA-SQ)\psi\\
&=(\partial+\dbar+A+Q)\psi-(-\partial+\dbar-A+Q)\psi\\
&=2(\partial+A)\psi.
\end{align*}
Therefore
\begin{align*}
(d*A)(X,JX)\psi
&=-2\dbar_X(A(X)\psi)+2A(X)^2\psi+2A(X)\partial_X\psi+2A(X)^2\psi\\
&=-2(\dbar_X(A(X)\psi)-A(X)\partial_X\psi)\\
&=-2(\dbar_XA)(X)\psi.
\end{align*}
\end{proof}

Now assume that $L$ is Willmore, and therefore $d*A=0$. This implies $\dbar
A=0$, and $A$ is holomorphic:
\begin{align*}
A\in H^0(K\End_-(H))=H^0(K\Hom_+(\bar{H},H)).
\end{align*}
As a consequence, see
Lemma~\ref{lemma:holohomo}, either $A\equiv0$, or the zeros of $A$ are
isolated, and there exists a line bundle $\tL\subset H$ such that
$\tL=\ker A$ away from the zeros of $A$. For local $\psi\in\Gamma(\tL)$ and
holomorphic $Y\in H^0(TM)$ we have
\begin{align*}
\underbrace{\dbar A}_{=0}(Y)\psi
=\dbar(\underbrace{A(Y)\psi}_{=0})-A(Y)\partial \psi.
\end{align*}
Therefore $\tL$ is invariant under $\partial$, like $L$ is invariant under
$\dbar$, see Remark~\ref{remark:holostrH}.
As above, we get a holomorphic structure on the complex {\em line} bundle
$K\Hom_+(\bar{H}/\tL,L)$ and $A$ defines a holomorphic section of this
bundle:
\begin{align*}
A\in H^0(K\Hom_+(\bar{H}/\tL,L)).
\end{align*}
\end{example}
\end{subsection}

\begin{subsection}{Spherical Willmore Surfaces}
We turn to the
\begin{proof}[Proof of Theorem~\ref{theorem:montiel}]
If $A\equiv0$ or $Q\equiv0$, then $L$ is a twistor projection by
Theo\-rem~\ref{theorem:superconformal}.

Otherwise we have the line bundle
$\tL$, and similarly a line bundle $\hL$ that coincides with the image of
$Q$ almost everywhere.

\begin{proposition}
We have the following holomorphic sections of complex holomorphic line
bundles:
\begin{alignat*}{2}
A&\in H^0(K\Hom_+(\bar{H}/\tL,L)),
\quad&Q&\in H^0(K\Hom_+(H/L,\bar{\hL})),\\
\delta_L&\in H^0(K\Hom_+(L,H/L)),
\quad&AQ&\in H^0(K^2\Hom_+(H/L,L))\\
&\text{ and if  $AQ=0$ then }
\quad&\delta_{\tL}&\in H^0(K\Hom_+(\bar{\tL},\bar{H}/\tL))
\end{alignat*}
\end{proposition}

We proved the statement about $A$. We give the (similar) proofs of the
others in the appendix.

The degree formula then
yields
\begin{align*}
\ord \delta_L&=\deg K-\deg L+\deg H/L \\
\ord (AQ)&=2\deg K-\deg H/L+\deg L\\
&=3\deg K-\ord\delta_L\\
&=6(g-1)-\ord\delta_L.
\end{align*}
For $M=S^2$, i.e. $g=0$, we get $\ord (AQ)<0$, whence $AQ=0$. Then
$\tL =\hL $, and
\begin{align*}
\ord A&=\deg K +\deg H/\tL +\deg L\\
\ord Q&=\deg K -\deg H/L-\deg \tL \\
\ord\delta_{\tL }&=\deg K +\deg\tL -\deg H/\tL.
\end{align*}
Addition yields
\begin{align*}
\ord\delta_{\tL }+\ord Q+\ord A
&=3\deg K-\deg H/L+\deg L\\
&=4\deg K-\ord\delta_L=-8-\ord\delta_L.
\end{align*}
It follows that $\ord\delta_{\tL }<0$, i.e. $\delta_{\tL }=0$,
and $\tL $ is $d$-stable, hence constant in $H=M\times\HH^2$. From
$AS=-SA=0$ we conclude $S\tL =\tL $. Therefore all mean
curvature spheres of $L$ pass through the fixed point $\tL $.
Choosing affine coordinates with $\tL =\infty$, all mean curvature
spheres are affine planes, and $L$ corresponds to a minimal surface in
$\RR^4$.
\end{proof}
\end{subsection}
\end{section}

\newpage

\begin{section}{Appendix}\label{subsection:supplement}
\begin{subsection}{The bundle $\tL$}
\begin{lemma}\label{lemma:Anull} If $L$ is is an immersed
holomorphic curve in $\HH P^1$  with $d*dS=0$ then
\begin{align*}
A|_L=0\iff A=0.
\end{align*}
\end{lemma}
\begin{proof}
 $0=d*dS=2d(A-Q)$ implies
\begin{align*}
dA=\frac{1}{2}d(A+Q)=Q\wedge Q+A\wedge A,
\end{align*}
see Lemma \ref{lemma:dAplusQ}. Since
$Q|_L=0$, the assumption $A|_L=0$ implies $dA|_L=0$. Then for
$\psi\in\Gamma(L)$
\begin{align*}
0=d(A\psi)=(dA)\psi-A\wedge d\psi=-A\wedge d\psi.
\end{align*}
Since $A|_L=0$, this implies
\begin{align*}
0=A\wedge\delta=A*\delta-*A\delta=-2SA\delta.
\end{align*}
But $L$ is an immersion. Therefore $A|_L=0=A\delta$ implies $A=0$.
The converse is obvious.
\end{proof}
\begin{lemma}\label{lemma:holohomo}
Given a holomorphic section $T\in H^0(\Hom(V,W))$, where $V,W$
are holomorphic complex vector bundles, there exist holomorphic
subbundles
\begin{align*}
V_0\subset V, \hat{W}\subset W
\end{align*}
such that $V_0=\ker
T$ and $\hat{W}=\im T$ away from a discrete subset.
\end{lemma}
\begin{proof}
Let $r:=\max \{\rank T_p\,|\,p\in M\}$ and $G:=\{p\,|\,\rank T_p=r\}$. This is
an open subset of
$M$. Let $p_0$ be a boundary point of $G$, an let $\psi_1,\ldots,\psi_n$ be
holomorphic sections of $V$ on a neighborhood $U$ of $p_0$. By a change of
indices we may assume that
$T\psi_1\wedge\ldots\wedge T\psi_r\not{\!\!\equiv} 0$. But this is a
holomorphic section of the  holomorphic bundle
$\Lambda^r W|_U$, and hence has isolated zeros, because $\dim_{\CC}M=1$.
We assume that $p_0$ is its only zero within $U$. Moreover,
there exist
$k\in\NN$, a holomorphic coordinate $z$  centered at $p_0$, and a
holomorphic section
$\sigma\in H^0(\Lambda^r W|_U)$ such that
\begin{align*}
T\psi_1\wedge\ldots\wedge T\psi_r=z^k\sigma.
\end{align*}
Off $p_0$ the section $\sigma$ is decomposable, and since the Grassmannian
$G_r(W)$ is closed in $\Lambda^r(W)$, it defines a section of $G_r(W)$,
i.e. an $r$-dimensional subbundle of $W|_U$ extending $\im T|_{U\backslash
p_0}$. The statement about the kernel follows easily using the fact that
$\ker T$ is the annihilator of $\im T^*:W^*\to V^*$.

\end{proof}

\begin{proposition}
Let $L$ be a (connected) Willmore surface in $\HH P^1$, and $A\not\equiv
0$. Then there exists a unique line bundle $\tL\subset H$ such that on an
open dense subset of $M$ we have:
\begin{align*}
\tL=\ker A \text{ and } H=L\oplus\tL.
\end{align*}
\end{proposition}
\begin{proof}
$A\in\Gamma(K\End_-(H))$ is a holomorphic section by
Example~\ref{example:Aholo}. By Lemma~\ref{lemma:holohomo} there exists
a line bundle $\tL$ such that $\tL=\ker A$ off a discrete set. Assume now that
$H|_U\neq L\oplus\tL$ on an open non-empty set $U\subset M$. Then $L=\tL$, and
$A|_L=0$ on $U$. But then
$A|_U=0$ by Lemma~\ref{lemma:Anull}. This is a contradiction, because the
zeros of
$A$ are isolated.
\end{proof}
\end{subsection}

\begin{subsection}{Holomorphicity and the Montiel theorem}

In this section $L$ denotes an immersed holomorphic curve in
$\HH P^1$.

\begin{remark}[Holomorphic Vector Fields]\label{remark:holoVF}
The tangent
bundle of a Riemann surface viewed as complex line bundle carries a holomorphic
structure:
\begin{align*}
\dbar_XY=\frac{1}{2}([X,Y]+J[JX,Y]).
\end{align*}
Note that this is tensorial in $X$. The vanishing of the Nijenhuis tensor
implies $\dbar J=0$. A vector field
$Y$ is called holomorphic if $\dbar Y=0$. This is equivalent with
$\dbar_YY=0=\dbar_{JY}Y$, but either of these conditions simply says
\begin{align*}
[Y,JY]=0.
\end{align*}
Any constant vector field in $\CC$ is therefore holomorphic, and a given
tangent vector to a Riemann surface can always be extended to a holomorphic
vector field.
\end{remark}

\begin{proposition}
Let $L$ be a Willmore surface in $\HH P^1$. We have the
following holomorphic sections of complex holomorphic line bundles:
\begin{alignat*}{2}
A&\in H^0(K\Hom_+(\bar{H}/\tL,L)),
\quad&Q&\in H^0(K\Hom_+(H/L,\bar{\hL})),\\
\delta_L&\in H^0(K\Hom_+(L,H/L)),
\quad&AQ&\in H^0(K^2\Hom_+(H/L,L)),\\
&\text{ and if  $AQ=0$ then }
\quad&\delta_{\tL}&\in H^0(K\Hom_+(\bar{\tL},\bar{H}/\tL)).
\end{alignat*}
\end{proposition}
For the proof we need

\begin{lemma}\label{lemma:curv}
 The curvature tensor of the connection
$\partial+\dbar$ on $H$ is given by
\begin{align}
R^{\partial+\dbar}=-(A\wedge A+Q\wedge Q),\label{eq:RAQ}
\end{align}
and for a holomorphic vector field $Z$ we have
\begin{align}
R^{\partial+\dbar}(Z,JZ)=2S(\dbar_Z\partial_Z-\partial_Z\dbar_Z).
\label{eq:Rpartial}
\end{align}
\end{lemma}
\begin{proof} In general, if $\nabla$ and $\tilde{\nabla}=\nabla+\omega$
are two connections, then
\begin{align*}
R^{\tilde{\nabla}}=R^\nabla+d^\nabla\omega+\omega\wedge\omega.
\end{align*}
We apply this to $\tilde\nabla=\partial+\dbar=d-(A+Q)$ and use
Lemma~\ref{lemma:dAplusQ}:
\begin{align*}
R^{\partial+\dbar}&=R^d-d(A+Q)+(A+Q)\wedge (A+Q)\\
&=-2(A\wedge A+Q\wedge Q)+(A\wedge A+Q\wedge Q)\\
&=-(A\wedge A+Q\wedge Q).
\end{align*}
Equation \eqref{eq:Rpartial} follows from
\begin{align*}
R^{\partial+\dbar}(Z,JZ)
&=(\partial_Z+\dbar_Z)(\partial_{JZ}+\dbar_{JZ})
 -(\partial_{JZ}+\dbar_{JZ})(\partial_Z+\dbar_Z)\\
&=S(\partial_Z+\dbar_Z)(\partial_{Z}-\dbar_{JZ})
 -S(\partial_{Z}-\dbar_{Z})(\partial_Z+\dbar_Z)\\
&=2S(-\partial_Z\dbar_{JZ}+\dbar_{Z}\partial_Z),
\end{align*}
because $\dbar_Z^2=0=\partial_Z^2$.
\end{proof}

\begin{proof}[Proof of the proposition]

The holomorphicity of $A$ was shown in example~\ref{example:Aholo}, and
that of
$Q$ can be shown in complete analogy.

$(H,S)$ is a holomorphic complex quaternionic vector bundle, and $L$ is a
holomorphic subbundle, see Remark~\ref{remark:holostrH}. Therefore $L$ and
$H/L$ are holomorphic complex quaternionic line bundles, and  the complex line
bundle $K\Hom_+(L,E/L)$ inherits a holomorphic structure. Then, for
local holomorphic sections $\psi$ in $L$ and $Z$ in $TM$,
\begin{align*}
(\dbar_Z\delta_L)(Z)\psi
&=\dbar_Z(\delta_L(Z)\psi)-\delta_L(\dbar_ZZ)\psi-\delta_L(Z)(\dbar_Z\psi)\\
&=\dbar_Z(\delta_L(Z)\psi)
=\dbar_Z(\pi_L d\psi(Z))\\
&=\pi_L \dbar_Z(d\psi(Z))
=\pi_L \dbar_Z(\partial_Z\psi).
\end{align*}
By \eqref{eq:RAQ} and \eqref{eq:Rpartial} we have
\begin{align*}
\dbar_Z\partial_Z\psi=\partial_Z\underbrace{\dbar_Z\psi}_{=0}
-\frac{1}{2}\underbrace{R^{\partial+\dbar}(Z,JZ)\psi}_{\in  L},
\end{align*}
hence
\begin{align*}
(\dbar_Z\delta_L)(Z)=0.
\end{align*}
Then also
\begin{align*}
(\dbar_{JZ}\delta_L)(Z)=S(\dbar_{Z}\delta_L)(Z)=0,
\end{align*}
and therefore $\dbar\delta_L=0$.

To prove the holomorphicity of $AQ\in\Gamma(K^2\Hom(H/L,\bar{\hL}))$, we first
note that
\begin{align*}
K^2\Hom(H/L,\bar{\hL})=\Hom_\CC(TM,\Hom_\CC(TM,\Hom_+(H/L,\bar{\hL})))
\end{align*}
carries a natural holomorphic structure. The rest follows from the
holomorphicity
of $A,Q$, and the product rule.

Finally we interpret $\delta_{\tL}$ as a section in
$K\Hom_+(\bar{\tL},\bar{H}/\tL)$. Note that the holomorphic structure on
$\bar{H}$ is given by $\partial$. From the holomorphicity of $A$ we find,
for $\phi\in\Gamma(\tL)$,

\begin{align*}
0=(\dbar A)\phi=\dbar(\underbrace{A\phi}_{=0})+A\partial\phi.
\end{align*}
This shows that $\tL$ is $\partial$-invariant. Moreover, it is obviously
invariant under $A$ and, as a consequence of $AQ=0$, also under $Q$. From
Lemma~\ref{lemma:curv} it follows that $\tL$ is invariant under
$R^{\partial+\dbar}$, and that for a local holomorphic vector field $Z$ and a
local holomorphic section $\phi$ of $\tL$,
\begin{align*}
\partial_Z\dbar_Z\phi=\dbar_Z\underbrace{\partial_Z\phi}_{=0}
+\frac{1}{2}\underbrace{SR^{\partial+\dbar}(Z,JZ)\psi}_{\in  \tL}.
\end{align*}
Then
\begin{align*}
(\dbar_Z \delta_{\tL})(Z)\phi
&=\partial(\delta_{\tL}(Z)\phi)-\delta_{\tL}(\dbar_Z
Z)\phi-\delta_{\tL}(Z)\partial_Z\phi\\
&=\partial_Z (\delta_{\tL}(Z)\phi)
=\partial_Z (\pi_{\tL}d\phi(Z))
=\pi_{\tL}\partial_Z (d\phi(Z))\\
&=\pi_{\tL}\partial_Z \dbar_Z\phi
=0.
\end{align*}
\end{proof}
\end{subsection}
\end{section}

\begin{section}{%
Epilogue}

In the presentation of the material given in this
course, I strictly focused on surfaces in $\HH P^1$,
though many concepts may also be considered for
surfaces in $\HH P^n$ or even for more general
situations. A significant difference in  higher codimensions is the lack
of a unique mean curvature sphere congruence as given by
Theorem~\ref{theorem:mcs}. As a consequence, the bundle $L$ will not carry a
natural holomorphic structure. But $L^{-1}$ will: see
Theorem~\ref{theorem:Linverse}.

The global theory of holomorphic sections and degree theory for complex
quaternionic line bundles is under construction, see e.g. \cite{icm}.
There one finds  a lower bound for the Willmore functional:
\begin{align*}
W(L)\ge -d+\ord\psi,
\end{align*}
for a nontrivial section $\psi\in H^0(L^{-1})$. Here $d:=\deg(L^{-1})$ is the
degree of $L^{-1}$. In \cite{BPP} a stronger inequality  will be shown under
certain non-degeneracy assumptions:
\begin{align*}
W(L)\ge h^0(h^0-d-1)
\end{align*}
where $h^0:=\dim H^0(L^{-1})$.

Another topic addressed in \cite{icm} is that of holomorphic structures
on paired complex quaternionic line bundles and  generalized
Weierstrass representations. Let
$L$ be a complex quaternionic line bundle with holomorphic structure
$D$. Then $KL^{-1}$ carries a unique holomorphic structure $\tilde{D}$ such
that, as quadratic forms,
\begin{align*}
d<\alpha,\psi>=<{\tilde{D}\alpha},J\psi>-<J\alpha, D\psi>
\end{align*}
for $\alpha\in \Gamma(KL^{-1}),\psi\in \Gamma(L)$. In this situation, the
Riemann-Roch Theorem,
\begin{align*}
\dim H^0(L)-\dim H^0(KL^{-1})=\deg(L)-g+1,
\end{align*}
 holds on compact Riemann surfaces.

Given holomorphic sections $\alpha\in H^0(KL^{-1}),\; \psi\in H^0(L)$ there
exists a local $f:M\to\HH$ such that $df=<\alpha,\psi>$ and $f$ is
conformal with right normal given by $J\psi=-\psi R$. Conversely,  any
conformal
$f$ can be obtained in this way: Put $L:=M\times\HH, J\psi:=-\psi R$,
and
$D1:=0$. Then $df=<\alpha,1>$  determines a holomorphic section of
$KL^{-1}$.

Besides Willmore surfaces, the family of isothermic surfaces fits perfectly
into
the present frame. Classically, in
$\RR^3$, they are defined by the property of carrying conformal curvature line
coordinates, or, equivalently, by the fact that their mean curvature
sphere congruence touches a second enveloping surface  conformally related
to, but with opposite orientation from the original one. Examples are Willmore
surfaces or constant mean curvature surfaces in 3-space. Isothermic surfaces
have been also defined in 4-space, see \cite{HJP}. In our setting, we call
 $f:M\to\HH P^1$ isothermic, if there exists a second surface
$g:M\to\HH P^1$ such that $df\wedge dg=0=dg\wedge df$. There is a quite
satisfactory generalization of the classical  Darboux transformation theory for
these surfaces, see \cite{BPP}.
\end{section}


\par\medskip\usebox{\address}

\begin{thebibliography}{}

\bibitem{bryant} Bryant, Robert.
{\em A duality theorem for Willmore surfaces.} J. Differential Geom. 20,
23-53 (1984)
%
\bibitem{ejiri} Ejiri, Norio.
{\em Willmore Surfaces with a Duality in $S^N$(1).} Proc. Lond. Math.
Soc., III Ser. 57, No.2, 383-416 (1988)
%
\bibitem{friedrich} Friedrich, Thomas.
{\em On Superminimal Surfaces.} Archivum math. 33, 41-56 (1997)
%
\bibitem{HJP} Hertrich-Jeromin, Udo; Pedit, Franz.
{\em Remarks on Darboux Transforms of Isothermic Surfaces.} Doc. Math. J.
2, 313 - 333 (1997).\\
(www.mathematik.uni-bielefeld.de/documenta/vol-02/vol-02.html)

\bibitem{kulkarni} Kulkarni, Ravi; Pinkall, Ulrich (Eds.).
{\em Conformal Geometry.} Vieweg, Braunschweig 1988
%
\bibitem{montiel} Montiel, Sebasti\'an.
{\em Spherial Willmore Surfaces in the Four-Sphere.} Preprint 1998
%
\bibitem{icm} Pedit, Franz; Pinkall, Ulrich.
{\em Quaternionic analysis on Riemann surfaces and differential
geometry.} Doc. Math. J. DMV, Extra Volume ICM 1998, Vol. II, 389-400.\\
(www.mathematik.uni-bielefeld.de/documenta/xvol-icm/05/05.html)
%
\bibitem{BPP} Pedit, Franz; Pinkall, Ulrich; et al.
{\em Quaternionic Holomorphic Geometry and Willmore surfaces in $S^4$.}
To appear. \\(Check
www-sfb288.math.tu-berlin.de/Publications/Preprints.html)
%
\bibitem{richter} Richter, J\"org.
{\em Conformal Maps of a Riemannian Surface into the Space of
Quaternions.} Dissertation, Berlin 1997
%
\bibitem{rigoli} Rigoli, Marco.
{\em The conformal Gauss map of Submanifolds of the Moebius Space.}
Ann. Global Anal. Geom 5, No.2, 97-116 (1987)
%
\bibitem{tenenblat} Tenenblatt, Keti.
{\em Transformations of Manifolds and Applications to Diffe\-rential
Equations.} Chapman\& Hall/CRC Press 1998
\end{thebibliography}
\end{document}